\documentclass[12pt,reqno]{amsart}
\usepackage{amscd, amssymb}

\DeclareMathOperator*{\tnull}{}

\swapnumbers
\newtheorem{thm}[subsection]{Theorem}
\newtheorem{cor}[subsection]{Corollary}
\newtheorem{lem}[subsection]{Lemma}
\newtheorem{prop}[subsection]{Proposition}

\theoremstyle{definition}
\newtheorem{defn}[subsection]{Definition}

\newcommand{\tr}{\sideset{^{t{}}}{}\tnull\nolimits\hskip-3.5pt}
\newcommand{\trA}{\sideset{^{t{}}}{}\tnull\nolimits\hskip-4.5pt A}
\renewcommand{\Sp}{\mathbf{Sp}(2n,\mathbb R)}
\newcommand{\quash}[1]{}
\newcommand{\parens}[1]{{\rm(}#1{\rm)}}
\renewcommand{\mod}[1]{\ (\text{\rm mod }#1)}

\numberwithin{equation}{subsection}

\begin{document}

\title{The moduli space of real abelian varieties with level structure}
\author{Mark Goresky${}^1$}\thanks{1.  School of Mathematics, Institute for Advanced Study,
Princeton N.J.  Research partially supported by NSF grants \# DMS 9626616 and DMS
9900324.}
\author{Yung sheng Tai${}^2$}\thanks{2.  Dept. of Mathematics, Haverford College, Haverford
PA.}
\begin{abstract}
 The moduli space of principally polarized abelian
varieties with real structure and with level $N=4m$ structure (with $m \ge 1$) is shown to
coincide with the set of real points of a quasi-projective algebraic variety defined over
$\mathbb Q$, and to
consist of finitely many copies of the quotient of the space $\mathbf{GL}(n, \mathbb
R)/\mathbf O(N)$ (of positive definite symmetric matrices) by the principal congruence
subgroup of level $N$ in $\mathbf{GL}(n, \mathbb Z).$
\end{abstract}\maketitle

\section{Introduction}
Let $\mathfrak h_n = \mathbf{Sp}(2n, \mathbb R)/\mathbf{U}(n)$ be the Siegel upperhalf space
of rank $n$.  The quotient space $\mathbf{Sp}(2n,\mathbb Z) \backslash\mathfrak h_n$ has three
remarkable properties:  (a)  it is the moduli space of principally polarized abelian
varieties, (b) it has the structure of a quasi-projective complex algebraic variety which is
defined over the rational numbers $\mathbb Q,$ and (c) it has a natural compactification (the
Baily-Borel Satake compactification) which is defined over the rational numbers. 

Now let $C_n = \mathbf{GL}(n,\mathbb R) /
\mathbf{O}(n)$ be the symmetric cone of positive definite symmetric matrices and let $Z =
\mathbf{GL}(n,\mathbb Z) \backslash C_n.$  One might ask whether similar statements hold for
$Z.$  Is $Z$ in some sense a moduli space for principally polarized real abelian varieties?
Does it admit the structure of a real algebraic variety, possibly defined over $\mathbb
Q$?   If so, does it admit a compactification which is also  defined over the rational
numbers?  The answer to all these questions is ``no''.  In fact,
Silhol \cite{Silhol} constructs the moduli space of real principally polarized abelian
varieties and he shows that it is a (topological) ramified covering of $Z$.  Moreover Silhol
constructs a compactification of this moduli space, analogous to the Baily-Borel
compactification.  However neither the moduli space nor this compactification has an algebraic
structure.  

In this paper we show that all three statements (a), (b), and (c)  above may be recovered if
we consider real abelian varieties with an appropriate level structure.  To be precise, let $N
= 4m$ for some $m \ge 1,$ let $\Gamma_{\ell}(N)$ be the principal congruence subgroup of
$\mathbf{GL}(n,\mathbb Z)$ of level $N$, and let $Z(N) = \Gamma_{\ell}(N)\backslash C_n.$
Then (a) the moduli space of real principally polarized abelian varieties with level $N$
structure consists of finitely many copies of $Z(N).$  These various copies are indexed by a
certain (nonabelian cohomology) set $H^1(\mathbb C/ \mathbb
R, \Gamma(N)).$  Morever, (b) this moduli space may be naturally identified with the real
points $X_{\mathbb R}$ of a quasi-projective algebraic variety $X$ which is defined over the
rational numbers.  See Theorems \ref{thm-real} and \ref{thm-moduli}.  Finally, (c) the variety
$X_{\mathbb C}$ has a compactification which is defined over $\mathbb Q.$

In fact the variety $X_{\mathbb C}$ is just the moduli space $\Gamma(N) \backslash \mathfrak
h_n$ of abelian varieties with principal polarization and level $N=4m$ structure.  (Here,
$\Gamma(N)$ denotes the principal congruence subgroup of $\mathbf{Sp}(2n,\mathbb Z)$ of level
$N.$)  The key observation is that there is an involution $\tau$ on $\mathbf{Sp}(2n,\mathbb
R)$ whose fixed point set is $\mathbf{GL}(n, \mathbb R),$ which passes to an anti-holomorphic
involution on $X_{\mathbb C}.$

Let $V=V_{\mathbb C}$ denote the Baily-Borel Satake compactification of $X_{\mathbb
C}$.  It is an algebraic variety defined over $\mathbb Q.$   Complex conjugation $\tau:
V_{\mathbb C} \to V_{\mathbb C}$ is an anti-holomorphic involution whose fixed point set is
the set of real points $V_{\mathbb R}.$  This set is compact and it contains
$X_{\mathbb R}$ as an open set.  One might ask whether
$V_{\mathbb R}$ is a compactification of $X_{\mathbb R}$, that is, whether $V_{\mathbb R} =
\overline{X}_{\mathbb R},$ where $\overline{X}_{\mathbb R}$ denotes the closure of $X_{\mathbb
R}$ in $V_{\mathbb C}.$  We have been able to prove this (\S \ref{sec-2k}) in the special case
$N=2^k$ (that is, for the principal congruence subgroup $\Gamma(2^k)$ of level $2^k$) for $k
\ge 2$, and we suspect it is otherwise false unless $n=1;$ cf \S \ref{sec-SL2}.  However we
have been able to show (for general $\Gamma(4m)$) that the complement $V_{\mathbb R} -
\overline{X}_{\mathbb R}$ has a high codimension. The variety $V_{\mathbb C}$ has a
stratification
\[ V_{\mathbb C} = V^0_{\mathbb C} \cup V^1_{\mathbb C} \cup\ldots \cup V^n_{\mathbb C}\]
such that each $V^r_{\mathbb C}$ is a union of arithmetic quotients $\Gamma \backslash
\mathfrak h_{n-r}$ of finitely many copies of the Siegel space of rank $n-r$.  The largest
stratum $V^0_{\mathbb C}$ is just $X_{\mathbb C}.$  Denote by $V^r_{\mathbb R}$ the
$\tau$-fixed points in $V^r_{\mathbb C}.$  In Proposition \ref{prop-BBcompact} we show that
\[ V^0_{\mathbb R} \cup V^1_{\mathbb R} \subset \overline{X}_{\mathbb R} \subset V_{\mathbb
R}. \] 

Although the moduli space $X_{\mathbb R}$ consists of finitely many disjoint copies of the
locally symmetric space $Z(N)$, the compactification $\overline{X}_\mathbb R$ is not a
disjoint union:  some of these copies of $Z(N)$ may become glued together along the boundary.

The authors are grateful to the Institute for Advanced Study in Princeton N.J. for its
hospitality and support.

\section{Symplectic Group}
\subsection{}
The symplectic group $G=\Sp$ may be realized as the group of $2n$ by $2n$ real
matrices
\begin{align*}
\Sp &= \Bigl\{ \Bigl( 
\begin{array}{cc} A & B \\ C & D \end{array}\Bigr) \ \Bigm| \  
\begin{array}{ll} \trA D - \tr C B = I \\
\trA C \text{,} \tr B D \text{ symmetric}
\end{array}\Bigr\} \\ 
&= \Bigl\{ \Bigl( 
\begin{array}{cc} A & B \\ C & D \end{array}\Bigr) \ \Bigm| \  
\begin{array}{ll}  A \tr D -  B \tr C = I \\
 A \tr B \text{,}  C \tr D \text{ symmetric}
\end{array}\Bigr\}
\end{align*}
If $J={\bigl( \begin{matrix} {0} &  {I_n} \\ {-I_n} &{0} \end{matrix} \bigr)}$ is the
standard symplectic form on $\mathbb R^{2n}$ then $g\in\Sp \iff \tr g  Jg=J.$
The inverse of such a symplectic matrix is
\begin{equation}\label{eqn-inverse}
\Bigl( \begin{array}{cc}
 A & B \\ C & D \end{array} \Bigr)^{-1} = 
\Bigl( \begin{array}{cccc}
\tr D  & -^tB \\ -^tC & \trA \end{array}\Bigr)
\end{equation}
Identify  $\mathbf {GL}(n,\mathbb R) \hookrightarrow \Sp$  with its image under the
embedding
\begin{equation*}
A \mapsto
\Bigl(
\begin{array}{cc} A & 0 \\ 0 & \trA^{-1} \end{array} \Bigr) 
\end{equation*}
A Cartan involution $\theta$ of $\Sp$ is given by $\theta g = Jg J^{-1}$, that is,
\begin{equation*}
\theta \Bigl( \begin{array}{cc} A & B \\ C & D \end{array} \Bigr) =
\Bigl( \begin{array}{cc} D & -C \\ -B & A \end{array} \Bigr)
\end{equation*}
Its fixed point set is the unitary group $K=\mathbf{U}(n)$ which
is embedded in the symplectic group by 
\begin{equation*}
A+iB \mapsto  
\Bigl( \begin{array}{cc} A & B \\ -B & A \end{array} \Bigr).
\end{equation*}
 
The symplectic group acts transitively on the Siegel upper halfspace 
\begin{equation*}\mathfrak h_n = \left\{Z=
X+iY \in M_{n\times n}(\mathbb C)\ \Bigm|\ \tr Z=Z \text{ and } Y>0  \right\}
\end{equation*} 
by fractional linear transformations:  if $g =\left( \begin{smallmatrix} 
A & B \\ C & D \end{smallmatrix} \right)$ then 
\begin{equation*}
gZ = (AZ+B)(CZ+D)^{-1}.\end{equation*} 
Then $\text{Stab}_G(iI) = \mathbf{U}(n),$ and $\mathfrak h_n \cong G/K$.  

\subsection{}\label{subsec-tau}
Let $I_{-} = \left( \begin{smallmatrix} -I & 0 \\ 0 & I \end{smallmatrix} \right) \in
\mathbf{GL}(2n, \mathbb R).$   Consider the involution $\tau:\Sp \to \Sp$ which is defined by
$\tau(g) = I_{-} g I_{-}$, that is, 
\begin{equation}\label{eqn-tau}
\tau\left( \begin{array}{cc} A & B \\ C & D \end{array} \right) = 
\left( \begin{array}{cc} A & -B\\ -C & D \end{array} \right)
\end{equation}
The following properties of $\tau$ are easily verified by direct calculation:
\begin{enumerate}
\item $\tau$ is an involution of $G$, i.e. $\tau(gh^{-1}) =
\tau(g)(\tau(h))^{-1}.$
\item  $\tau(g)=g \iff g \in \mathbf{GL}(n,\mathbb R)  $
\item $\tau\theta = \theta\tau$ so $\tau (\mathbf U(n)) = \mathbf U(n)$
\item If $C=A+iB\in \mathbf U(n)$ then $\tau(C) = \bar C= A-iB.$
\end{enumerate}

It follows that $\tau$ passes to an involution (which we denote by the same letter)
$\tau:\mathfrak h_n \to \mathfrak h_n$ on the Siegel space such that
\[ \tau(g\cdot Z) = \tau(g)\cdot \tau(Z)\] for all $g\in \Sp$ and for all $Z \in \mathfrak
h_n.$  In fact $\tau:\mathfrak h_n \to \mathfrak h_n$ is the antiholomorphic involution given
by
\begin{equation*} \tau (Z) = -\overline{Z}. \end{equation*}
Its fixed point set is the orbit
\begin{equation}\label{eqn-iCn}iC_n= {\mathbf GL}(n,\mathbb R).iI\subset M_{n\times n}(\mathbb
C)\end{equation}
of the general linear group.  Here, $C_n \cong {\mathbf GL}(n,\mathbb R)/\mathbf O(n)$ denotes
the cone of positive definite symmetric real matrices.
\subsection{Notation}  For $\gamma \in G$ and $Z \in \mathfrak h_n$ we will usually write
$\tilde \gamma = \tau(\gamma)$ and $\tilde Z = \tau(Z).$  For any $\gamma \in \Sp$, define the
locus of $\gamma$-real points (or, more precisely, the
set of $(\gamma,\tau)$-real points),
\begin{equation}\label{eqn-hn}
\mathfrak h_n^{\gamma} = \{ Z \in \mathfrak h_n |\ \gamma Z = \tau(Z) \} \end{equation}
If $\Gamma \subset \Sp$ is an arithmetic subgroup such that $\tau\Gamma = \Gamma$, define the
set $\mathfrak h_n^{\Gamma}$ of $\Gamma$-real points to be the union $\bigcup_{\gamma
\in \Gamma} \mathfrak h_n^{\gamma}.$

\section{Galois cohomology}\label{sec-Galois}
\subsection{}  The involution $\tau$ of (\ref{eqn-tau}) may be considered as an action of
$\text{Gal}(\mathbb C/ \mathbb R) = \{ 1, \tau \}$ on $\text{\bf Sp}(2n, \mathbb R)$.  Its
restriction to the unitary group $K=\mathbf U(n) \subset \Sp$ coincides with complex
conjugation.  For any $\gamma \in \Sp$ let $f_{\gamma}:\text{Gal}(\mathbb C/ \mathbb R) \to
\Sp$ be the mapping $f_{\gamma}(1)=I$ and $f_{\gamma}(\tau) = \gamma.$  Then $f_{\gamma}$ is a
1-cocycle iff $\gamma \tilde\gamma = I$, and $f_{\gamma}$ is a coboundary iff there exists $h
\in \Sp$ so that $\gamma = \tilde hh^{-1}.$  In this section we will prove the following:

\begin{prop}\label{prop-Galois}  With respect to the above action \parens{\ref{eqn-tau}}, the
(nonabelian) Galois cohomology sets $H^1(\mathbb C/\mathbb R, \mathbf{U}(n))$ and $H^1
(\mathbb C/\mathbb R, \Sp)$ are trivial.  \end{prop}

\begin{lem}\label{lem-squares}
If $g \in \mathbf U(n)$ and $g \tilde g = I$  then there exists
$\delta \in \mathbf U(n)$ such that $g = \tilde{\delta}\delta^{-1}.$
\end{lem}

\subsection{Proof}  For completeness we give a short proof of this well known result, which is
the first part of Proposition \ref{prop-Galois}.
Set $g = \bigl( 
\begin{smallmatrix} A & B \\ -B & A \end{smallmatrix}\bigr).$  Calculating $g \tilde g =I$
gives $A^2 + B^2 = I$ and $AB = BA.$  Comparing with (\ref{eqn-inverse}) gives $\trA = A$ and
$\tr B = B.$ 
Such commuting symmetric matrices may be simultaneously diagonalized by an orthogonal matrix
$h\in \mathbf O(n).$ In other words, $h(A+iB)h^{-1} = \Lambda\in \mathbf U(n)$ is a diagonal
matrix.  Choose a square root, $\mu = \sqrt{\bar\Lambda} \in \mathbf U(n)$ (by choosing a
square root of each diagonal entry) and set $\delta = h^{-1}\mu h.$ 
Since $\tilde h = h$ and $\bar \mu = \mu^{-1}$, we find
\begin{equation*}
\tilde\delta\delta^{-1} = h^{-1}\tau(\mu)\mu^{-1} h = h^{-1}\bar\mu \bar\mu h = h^{-1}\Lambda
h = A+iB.  \qed
\end{equation*}

\begin{lem}\label{lem-coboundary}  Let $\gamma \in \Sp$ and suppose that $f_{\gamma}$ is a
1-cocycle \parens{i.e., $\gamma \tilde\gamma = I$}.  Then $f_{\gamma}$ is a coboundary iff
$\mathfrak h_n^{\gamma} \ne \phi.$ \end{lem}

\subsection{Proof}
Suppose $f_{\gamma}$ is a coboundary, say $\gamma = \tilde h h^{-1}.$  Then $\gamma h iY =
\tau(hiY)$ for any $Y \in C_n$, so $\mathfrak h_n^{\gamma} \ne \phi.$  On the other hand,
suppose $Z \in \mathfrak h_n^{\gamma}.$  Write $Z = h(iI)$ for some $h \in \Sp.$  Then 
$\tilde h^{-1} \gamma h(iI) = \tilde h^{-1}\gamma Z = \tilde h^{-1} \tilde Z = \tau(iI) = iI$
hence $\alpha = \tilde h^{-1} \gamma h \in \mathbf U(n).$  Moreover, $\alpha \tilde \alpha =
\tilde h^{-1} \gamma \tilde \gamma \tilde h = I$ so by Lemma \ref{lem-squares} there exists
$\delta \in \mathbf U(n)$ with $\alpha = \tilde \delta \delta^{-1}.$  Therefore $\gamma =
\tau(h \delta) (h \delta)^{-1}.$  \qed

\subsection{} Suppose that $\gamma = \left( \begin{smallmatrix} A & B \\ C & D
\end{smallmatrix} \right) \in \Sp$ and that $\gamma \tilde \gamma = I.$  Then
\begin{equation}\label{eqn-conditions}
 A = \tr{D},\ C = \tr{C},\ B = \tr{B},\
A^2 - BC = I,\ CA = \tr{A}C, \text{ and } AB = B \tr{A}. \end{equation}

\begin{lem}\label{lem-diag}  Suppose  that $\gamma = \left( \begin{smallmatrix} 0 & -C^{-1} \\
C & 0 \end{smallmatrix} \right) \in \Sp$  and that $C$ is diagonal.  Then the cocycle
$f_{\gamma}$ is cohomologically trivial. \end{lem}

\subsection{Proof} 
By Lemma \ref{lem-coboundary} it suffices to show that there exists $Z \in
\mathfrak h_n$ such that $\gamma \cdot Z = - C^{-1} Z^{-1} C^{-1} = - \overline{Z}.$
If $C = \text{diag}(c_1,c_2,\ldots, c_n)$ then a solution is given by
$Z = iY$ where $Y=\text{diag}(y_1,y_2,\ldots,y_n)$ and $y_i = |c_i|^{-1}.$  \qed

\begin{lem}\label{lem-zero}
  Suppose $\gamma = \left( \begin{smallmatrix} A & B \\ C & D \end{smallmatrix}
\right) \in \Sp.$  Suppose that $\gamma \tilde \gamma =I$ and that $C=0.$  Then
the 1-cocycle $f_{\gamma}$ is cohomologically trivial.\end{lem}

\subsection{Proof}
The cocycle $f_{\gamma}$ is equivalent to $f_{\gamma'}$ where
\begin{equation*}
\gamma ' = \left( \begin{matrix} A & 0 \\ 0 & \tr{A} \end{matrix} \right) = \tilde h \gamma
h^{-1} \end{equation*}
and
\begin{equation*}
h = \left( \begin{matrix} I & \frac{1}{2} A^{-1}B \\ 0 & I \end{matrix} \right)
\end{equation*}
using (\ref{eqn-conditions}).  Then $A^2 = I$ so $A$ is $\mathbf{GL}(n, \mathbb R)$-conjugate
to a diagonal matrix $\Lambda$ (with eigenvalues $\pm 1$).  Therefore $f_{\gamma'}$
is equivalent to the 1-cocycle defined by the element $\left( \begin{smallmatrix} \Lambda & 0
\\ 0 & \Lambda^{-1} \end{smallmatrix} \right) \in \mathbf U(n).$  It follows from Lemma
\ref{lem-squares} that $f_{\gamma'}$ is cohomologically trivial.  \qed

\subsection{Proof of Proposition \ref{prop-Galois}}
Let $\gamma' = \left(\begin{smallmatrix}A' & B' \\C'&D'\end{smallmatrix}\right)$ and suppose
that $f_{\gamma'}$ is a cocycle.  By (\ref{eqn-conditions}), $C'$ is symmetric so there exists
$k \in \mathbf O(n)$ so that $C = k C' \tr{k}= \left(\begin{smallmatrix} C_1 & 0 \\ 0 & 0
\end{smallmatrix} \right)$ where $C_1 \in M_{r\times r}(\mathbb R)$ is
diagonal and nonsingular.  Let $h_1 = \left( \begin{smallmatrix} k & 0 \\ 0 & k
\end{smallmatrix} \right)$.  Then $f_{\gamma'}$ is cohomologous to $f_{\gamma}$ where 
\begin{equation*}
\gamma= \tilde h_1\gamma' h_1^{-1} = \left( \begin{matrix} A & B \\ C & D \end{matrix}
\right) \end{equation*}
and (\ref{eqn-conditions}) holds.  Write $A = \left( \begin{smallmatrix} A_1 & A_2 \\ A_3 &
A_4 \end{smallmatrix} \right)$ where $A_1\in M_{r\times r}(\mathbb R)$ and $A_4 \in
M_{(n-r)\times (n-r)}(\mathbb R),$ and similarly for $B$ and $D.$  Then $CA = \tr{A}C$ implies
that $A_2 = A_3 = 0.$

Set $h_2 = \left( \begin{smallmatrix} I & x \\ 0 & I \end{smallmatrix} \right)$ where 
\begin{equation*}
x = \left( \begin{matrix} C_1^{-1} \tr{A_1} & 0 \\ 0 & 0 \end{matrix} \right). 
\end{equation*}
Then $A^2 - BC = I$ implies that $f_{\gamma}$ is cohomologous to $f_g$ where
\begin{equation*}
g=\tilde h_2 \gamma h_2^{-1} = \left( \begin{array}{c|c}
A-xC & -Ax +B +xCx -x\tr{A} \\ \hline  C & -Cx + \tr{A}
\end{array} \right)
= \left( \begin{array}{cc|cc} 0 & 0 & -C_1^{-1} & 0 \\
0 & A_4 & 0 & B_4 \\ \hline
C_1 & 0 & 0 & 0 \\
0 & 0 & 0 & \tr{A_4} \end{array} \right).
\end{equation*}
Since $g \tilde g = I$ we may apply Lemma \ref{lem-diag} and Lemma \ref{lem-zero} to find $Z_r
\in \mathfrak h_r$ and $Z_{n-r} \in \mathfrak h_{n-r}$ so that 
\begin{equation*}
\left(\begin{matrix} 0 & -C_1^{-1} \\
C_1 & 0 \end{matrix} \right) Z_r = - \overline{Z}_r \text{ and }
\left(\begin{matrix} A_4 & B_4 \\ 0 & \tr{A_4} \end{matrix}\right)
Z_{n-r} = - \overline{Z}_{n-r}.\end{equation*}  Set $Z = \left( \begin{smallmatrix} Z_r & 0 \\
0 & Z_{n-r} \end{smallmatrix} \right).$  Then $gZ = - \overline{Z}$ hence the
corresponding cocycle is cohomologically trivial.  \qed


\begin{prop}\label{lem-Fgamma}
Let $\gamma \in \mathbf{Sp}(2n, \mathbb R)$ and suppose that $\tilde\gamma \gamma$ is
contained in a torsion-free subgroup of $\mathbf{Sp}(2n,\mathbb R).$  Then the following three
statements are equivalent:  \begin{enumerate}
\item $\mathfrak h_n^{\gamma} \ne \phi$ 
\item $f_{\gamma}$ is a cocycle (i.e.,  $\tilde \gamma \gamma = I $).
\item There exists $g \in \Sp$ such that $\gamma = \tilde
gg^{-1}.$ 
\end{enumerate}
In this case $\mathfrak h_n^{\gamma} = giC_n.$  Moreover such an element $g$ is
uniquely determined up to multiplication on the right by an element of $\mathbf{GL}(n,\mathbb
R).$   If $\Gamma \subset \Sp$ is a torsion free arithmetic group and if $\gamma_1, \gamma_2
\in\Gamma$ are distinct then $\mathfrak h_n^{\gamma_1} \cap \mathfrak h_n^{\gamma_2} = \phi.$
\end{prop}
\subsection{Proof}
Suppose $Z \in \mathfrak h_n^{\gamma} \ne \phi.$  Then $\tilde\gamma \gamma Z = \tilde \gamma
\widetilde{Z} = Z$ so $\tilde\gamma \gamma$ is torsion.  By hypothesis this implies that
$\tilde\gamma \gamma = I$, so (1) implies (2).  If  $f_{\gamma}$ is a 1-cocycle then by
Proposition \ref{prop-Galois} it is a coboundary.  Then Lemma
\ref{lem-coboundary} provides $g \in \Sp$ so that $\gamma = \tilde g g^{-1},$ which shows that
(2) implies (3).  It is easy to see that (3) implies (1).

If $\gamma = \tilde g g^{-1} = \tilde g'(g')^{-1}$ then $\tau$ fixes $\alpha=(g')^{-1}g.$
Hence $\alpha \in \mathbf{GL}(n,\mathbb R)$, which proves the next statement.

Finally, suppose $\gamma_1,\gamma_2\in\Gamma$ which is torsion-free.  If there exists $Z \in
\mathfrak h_n^{\gamma_1} \cap \mathfrak h_n^{\gamma_2}$ then $\gamma_1Z = \tau(Z) = \gamma_2Z$
so $\gamma_2^{-1}\gamma_1$ fixes $Z$, hence $\gamma_1 = \gamma_2.$ \qed

\section{The variety $\Gamma(4m)\backslash \mathfrak h_n$}
\subsection{}  In this section we will establish a strengthening of Proposition
\ref{lem-Fgamma} which applies to (certain) arithmetic subgroups of $\Sp.$
We shall be concerned with the following arithmetic groups: \begin{itemize}
\item $\Gamma(1) = \mathbf{Sp}(2n, \mathbb Z)$
\item $\Gamma(N) = \{ \gamma \in \Gamma(1) |\ \gamma \equiv I (\text{mod } N) \}$
\item $\Gamma_{2m}(2) = 
\Bigl\{\Bigl(\begin{array}{cc}  A & B \\ C & D \end{array} \Bigr) \in \Gamma \Bigm|
\begin{array}{cc}
A,D \equiv I\ (\text{mod }2) \\
B,C \equiv 0\ (\text{mod }2m) \end{array}\Bigr\}$
\item $\Gamma_{\ell}(1) = \mathbf{GL}(n,\mathbb Z)$
\item $\Gamma_{\ell}(N) = \{ \gamma \in \Gamma_{\ell}(1) |\ \gamma \equiv I 
(\text{mod } N) \}$
\end{itemize}
(The analogous construction of $\Gamma_m(2)$ does not yield a group unless $m$ is even.)
Throughout the remainder of the paper we will be interested in the case $N=4m$ and $m \ge 1.$ 
The following lemma (which depends on Proposition \ref{lem-surjective} below) is not essential
for this paper however it describes the relations between the above groups.

\begin{lem}\label{lem-doublecoset}  The short exact sequence
\begin{equation}\label{eqn-shortexact}
1 \to \Gamma(4m) \to \Gamma_{2m}(2) \to \Gamma(4m) \backslash \Gamma_{2m}(2) \to 1 
\end{equation}
induces a bijection
\begin{equation*}
H^1(\mathbb C / \mathbb R, \Gamma(4m)) \cong \Gamma(4m) \backslash \Gamma_{2m}(2) /
\Gamma_{\ell}(2).  \end{equation*}\end{lem}

\subsection{Proof}  We claim that $H^1(\mathbb C/ \mathbb R, \Gamma(4m)) \to H^1(\mathbb C/
\mathbb R, \Gamma_{2m}(2))$ is trivial.  For if $\gamma \in \Gamma(4m)$ and if $f_{\gamma}$ is
a 1-cocycle then by Proposition \ref{lem-Fgamma} (below), $\mathfrak h_n^{\gamma} \ne \phi.$
Then Proposition \ref{lem-surjective} implies that $\gamma = \tilde g g^{-1}$ for some $g \in
\Gamma_{2m}(2)$, which proves the claim.  Next, observe that $\text{Gal}(\mathbb C/ \mathbb
R)$ acts trivially on the quotient group $\Gamma(4m)\backslash\Gamma_{2m}(2).$  (For if $g \in
\Gamma_{2m}(2)$ then $g - \tilde g \equiv 0 \mod{4m}$ so $ \tilde g g^{-1} \equiv I
\mod{4m}$.)   Therefore the long exact cohomology sequence associated to
(\ref{eqn-shortexact}) is:
\begin{equation*}
\begin{CD}
H^0(\Gamma(4m)) &\to& H^0(\Gamma_{2m}(2)) &\to&
H^0(\Gamma(4m)\backslash \Gamma_{2m}(2)) &\to& H^1(\Gamma(4m)) &\to& 1 \\
|| && || && || \\
\Gamma_{\ell}(4m) && \Gamma_{\ell}(2) && \Gamma(4m)\backslash \Gamma_{2m}(2)
\end{CD}\end{equation*}
which completes the proof of Lemma \ref{lem-doublecoset}. \qed

Let $X= \Gamma(N) \backslash \mathfrak h_n.$  Then $X$ has the structure of a quasi-projective
complex algebraic variety.  The involution $\tau(Z) = - \overline{Z}$ passes to an
anti-holomorphic involution on $X$ and defines a real structure (let us call it the
$\tau$-real structure) on $X.$  In this section we will determine the real points of $X.$
(It is a theorem of Shimura \cite{Shimura} that the Baily-Borel compactification $V$ of $X$
even admits the structure of an algebraic variety defined over the rational numbers.  See \S
\ref{sec-BB}.)  Define
\[ S = \bigcup_{\gamma \in \Gamma(4m)} \mathfrak h_n^{\gamma} \subset \mathfrak h_n\]
to be the set of all $\Gamma(4m)$-real points in $\mathfrak h_n.$  Denote by
$\Gamma(4m)\backslash S$ the image of $S$ in $X = \Gamma(4m) \backslash \mathfrak h_n.$  For
$g \in \Gamma_{2m}(2)$ set ${}^g\Gamma_{\ell}(4m) = g\Gamma_{\ell}(4m)g^{-1}.$  In this
section we will prove the following
\begin{thm}\label{thm-real}
The set $X_{\mathbb R}$ of real points of $X$ is precisely $\Gamma(4m) \backslash S.$  It
consists of the disjoint union
\begin{equation}\label{eqn-copies}
X_{\mathbb R}= \coprod_g {}^g\Gamma_{\ell}(4m)\backslash giC_n \end{equation}
of finitely many copies of $\Gamma_{\ell}(4m) \backslash C_n,$ 
indexed by elements
\begin{equation}\label{eqn-doublecosets}
g\in\Gamma(4m)\backslash\Gamma_{2m}(2)/\Gamma_{\ell}(2) = H^1(\mathbb C/ \mathbb
R, \Gamma(4m)).\end{equation}
The copy indexed by $g \in \Gamma_{2m}(2)$ corresponds to the cohomology class of $f_{\gamma}$
where $\gamma = \tilde g g^{-1} \in \Gamma(4m).$
\end{thm}
In order to identify the set of real points $X_{\mathbb R}$ we will make repeated
use of the following observation:  {\it the image $[Z]\in X$ of a point $Z\in \mathfrak h_n$
lies in $X_{\mathbb R}$ iff  there exists $\gamma \in \Gamma(4m)$  such that $\gamma Z = -\bar
Z$} (i.e. such that $Z \in \mathfrak h_n^{\gamma}$).  The following lemma will be used
throughout the next few sections, and it also accounts for the indexing set in
(\ref{eqn-doublecosets}).

\begin{lem}\label{lem-divisibility}
Let $g\in \mathbf{Sp}(2n,\mathbb Z).$ Fix $m \ge 1.$\begin{enumerate}
\item If $g \in \Gamma(m)$ then $\tilde g g^{-1} \in \Gamma(2m).$  
\item If $g \in \Gamma(2m)$ then $\tilde g g \in \Gamma(4m).$ 
\item If $\tilde gg^{-1} \in \Gamma(4m)$ then $ g = \beta u$
for some $\beta \in \Gamma_{2m}(2)$ and for some $u \in \mathbf{GL}(n,\mathbb Z).$
\end{enumerate}\end{lem}

\subsection{Proof}  Let $g = \left( \begin{smallmatrix} A&B \\ C& D \end{smallmatrix}\right).$
Since $\tilde g = g - \bigl( \begin{smallmatrix} 0 & 2B \\ 2C & 0 \end{smallmatrix} \bigr)$
we have
\begin{equation}\label{eqn-trick}
 \tilde g g^{ -1} = I - \Bigl( \begin{array}{cc} -2B\tr C & 2B\trA \\ 2C\tr D & -2C\trA
\end{array}\Bigr).
\end{equation}
This gives part (1). If $g = I + 2mg'$ then $\tilde g + g$ is even, and $\tilde g g = I +
2m(\tilde g' + g') + 4m^2 \tilde g' g$ which proves part(2).  
Now suppose that $\tilde g g^{-1} \in \Gamma (4m)$ so that $B\tr C \equiv B \trA \equiv C \tr
D\equiv 0 (\text{mod } 2m).$  Then 
\begin{align*}
B &= B(\trA D - \tr C B) = (B\trA)D - (B \tr C) B \equiv 0\ (\text{mod } 2m)\\
\tr C &= (\trA D - \tr C B)\tr C = \trA (D\tr C) - \tr C (B \tr C) \equiv 0\ (\text{mod }
2m). 
\end{align*}
Moreover, since $C\tr B \equiv 0 (\text{mod }(2m))$ and $A\tr D - B \tr C = I$ we see that
$A\tr D \equiv I (\text{mod } 2)$.  Hence $A$ is invertible $(\text{mod } 2).$  But reduction
$(\text{mod } 2)$ is a surjective mapping $\mathbf{GL}(n,\mathbb Z) \to \mathbf{GL}(n, \mathbb
Z/(2))$ hence there exists $U \in \mathbf{GL}(n,\mathbb Z)$ so that $U \equiv A (\text{mod }
2)$ from which it also follows that $U^{-1}\equiv \tr D (\text{mod } 2)$.  Let $u=\bigl(
\begin{smallmatrix}U&0\\0&\tr U^{-1}
\end{smallmatrix}\bigr)\in \Sp.$Then
\begin{equation*}
\beta=gu^{-1}= \Bigl(\begin{matrix} AU^{-1} &B\tr U \\ CU^{-1} & D \tr U \end{matrix}\Bigr) 
\end{equation*} 
is in $\Gamma_{2m}(2)$ as claimed.  Conversely, if $\beta\in \Gamma_{2m}(2)$ and $u\in
\mathbf{GL}(n,\mathbb Z)$ then $\tau(\beta u) (\beta u)^{-1} = \tilde\beta\beta^{-1}$ which is
easily seen (using \ref{eqn-trick}) to lie in $\Gamma(4m).$  \qed

The following proposition is a consequence of the theorem of Silhol \cite{Silhol} (Thm.~1.4)
and Comessatti, and will be proven in \S \ref{sec-Comessatti}.
\begin{prop}\label{add}  Let $\gamma \in \Gamma(2)$ and suppose that $Z \in \mathfrak h_n$ is
not fixed by any element of $\mathbf{Sp}(2n,\mathbb Z)$ other than $\pm I.$  Suppose that
$\tilde Z=\gamma Z.$  Then there exists $h \in \mathbf{Sp}(2n,\mathbb Z)$ such that $\gamma =
\tilde h h^{-1},$ hence $\mathfrak h_n^{\gamma} = h \cdot iC_n.$ \end{prop}

\begin{cor}\label{lem-surjective} Suppose $\gamma \in \Gamma(4m)$ and $\mathfrak h_n^{\gamma}
\ne \phi.$  Then there exists $g\in \Gamma_{2m}(2)$ such that $\gamma = \tilde gg^{-1}$ (hence
$\mathfrak h_n^{\gamma} =giC_n$).
\end{cor}

\subsection{Proof}
Let $Z \in \mathfrak h_n^{\gamma}$ be a point which is not fixed by any $g\in
\mathbf{Sp}(2n,\mathbb Z)$ other than $\pm I.$  Such points exist, and even form a dense
subset of $\mathfrak h_n^{\gamma}$ because by Proposition \ref{lem-Fgamma}, the set $\mathfrak
h_n^{\gamma}$ is a translate (by some element of $\Sp$) of the cone $iC_n.$ By Proposition
\ref{add}, there exists $h \in \mathbf{Sp}(2n,\mathbb Z)$ so that $\gamma = \tilde hh^{-1}.$
By Lemma \ref{lem-divisibility} this implies that $h = \beta u$ for some $\beta \in
\Gamma_{2m}(2)$ and some $u \in \mathbf{GL}(n,\mathbb Z).$  Hence $\gamma = \tilde hh^{-1} =
\tilde\beta\beta^{-1}$ as claimed.  \qed

\begin{lem}\label{lem-cosets}
Let $g_1,g_2\in \Gamma_{2m}(2).$  Suppose there exist $Z_1\in g_1iC_n$, $Z_2\in g_2iC_n$, and
$\gamma \in \Gamma(4m)$ such that $\gamma Z_1 = Z_2.$  Then there exist $\alpha \in
\Gamma_{\ell}(2)$ such that $g_2 = \gamma g_1 \alpha.$
\end{lem}

\subsection{Proof}
By section \ref{subsec-tau} (3) we have a commutative diagram
\begin{equation*}
\begin{CD}
Z_1 @>{\gamma}>> Z_2 \\
@V{\tilde g_1g_1^{-1}}VV @VV{\tilde g_2g_2^{-1}}V \\
\tilde Z_1 @>>{\tilde\gamma}> \tilde Z_2
\end{CD}\end{equation*}
so the element $\theta = g_1\tilde g_1^{-1}\tilde\gamma^{-1}\tilde g_2g_2^{-1}\gamma$ fixes
$Z_1$.  By lemma \ref{lem-divisibility}, $\theta \in \Gamma(4m)$ which is torsion-free.  So
$\theta = I,$ which implies
\begin{equation*}
\gamma = g_2\tilde g_2^{-1}\tilde\gamma \tilde g_1g_1^{-1}
\end{equation*}
Therefore, 
$\tau(g_2^{-1}\gamma g_1) =g_2^{-1}\gamma g_1 $
which implies that $g_2^{-1}\gamma g_1 \in \mathbf{GL}_n(\mathbb R) \cap \Gamma_{2m}(2) =
\Gamma_{\ell}(2).$  \qed

\subsection{Proof of Theorem \ref{thm-real}}
For $\gamma \in \Gamma(4m)$ define $X^{\gamma} \subset X = \Gamma(4m) \backslash \mathfrak
h_n$ to be the image of $\mathfrak h_n^{\gamma}$ in $X$.  For each $g\in \Gamma_{2m}(2)$ the
corresponding element $\gamma = \tau(g)g^{-1}$ lies in $\Gamma(4m)$, and $\mathfrak
h_n^{\gamma} = giC_n.$  (This follows from lemma \ref{lem-divisibility} and Proposition
\ref{lem-Fgamma}.)  The resulting mapping
\begin{equation*}
\bigcup_{g\in\Gamma_{2m}(2)} giC_n \longrightarrow{}{} \bigcup_{\gamma \in
\Gamma(4m)}\mathfrak h_n^{\gamma}
\end{equation*}
is surjective by Corollary \ref{lem-surjective}.  For each $g\in \Gamma_{2m}(2)$ the
composition $giC_n \to \mathfrak h_n^{\gamma} \to X^{\gamma}$ is surjective and induces an
isomorphism ${}^g\Gamma_{\ell}(4m) \backslash giC_n \to X^{\gamma}$ (where
${}^g\Gamma_{\ell}(4m) = g\Gamma_{\ell}(4m)g^{-1}.$)  By lemma \ref{lem-cosets} this
determines a bijection
\begin{equation*}
\coprod_{g\in\Gamma(4m)\backslash\Gamma_{2m}(2)/\Gamma_{\ell}(2)}
 {}^g\Gamma_{\ell}(4m)\backslash giC_n \longrightarrow
\Gamma(4m) \backslash \bigcup_{\gamma \in \Gamma(4m)}\mathfrak h_n^{\gamma} = X_{\mathbb R}
\end{equation*}

\section{Upper half-plane}\label{sec-SL2}  For the upper half-plane $\mathfrak h_1$ it is 
possble to do a little better since the group $\Gamma(2)/\pm I$ acts freely.
\begin{lem}
A point $Z \in \mathfrak h_1$ is $(\Gamma(2),\tau)$-real if and only if $gZ$ is
$(\Gamma(2),\tau)$-real, for every $g\in SL(2,\mathbb Z).$
\end{lem}
\subsection{Proof.}  Suppose $\tau(Z) = \gamma Z$ for some $\gamma \in \Gamma(2).$  Then
\begin{equation*}
\tau(gz) = \tau(g)\tau(z) = (\tau(g)g^{-1})(g \gamma g^{-1}) gZ \in \Gamma(2)\cdot \Gamma(2)
\cdot gZ.  \qed
\end{equation*}
\begin{lem} \label{lem-upperhalfplane}
 If a point $y \in \mathfrak h_1$ is $(\Gamma(2),\tau)$-real, then 
there exists $\beta \in SL(2,\mathbb Z)$ and there exists $Y \in C_1 = 
\mathbb R_{+}$ so that $y = \beta\cdot iY.$ \end{lem}

\subsection{Proof.}  Say $\widetilde y = \gamma y$ for some $\gamma \in \Gamma(2).$ We
consider two cases:  (a) when $y$ is not fixed by any element of  $SL_2(\mathbb Z)$  other
than $\pm I$ and (b) when $y$ is fixed by some other nontrivial element of $SL_2(\mathbb Z).$ 

In case (a), the result is just Proposition \ref{add}.
Now consider case (b).  There are two classes of elements which are fixed
by nontrivial subgroups of $SL_2(\mathbb Z)$.  These are the translates (by elements of
$SL_2(\mathbb Z)$) of $i$ and the translates of $\omega=e^{2\pi i/3}.$  In the first case, $y
=\beta\cdot i$ (for some $\beta\in SL_2(\mathbb Z)$) and the conclusion follows.  In the
second case, the stabilizer of $\omega$ is the subgroup  
\begin{equation*} S= \left\{ \pm I, \pm \left(\begin{array}{rr}
0 & -1\\1 & 1 \end{array}\right), \pm \left( \begin{array}{rr}
-1 & -1\\ 1 & 0 \end{array}\right) \right\}.\end{equation*}
Suppose $y = g\omega$ (for some $g\in SL_2(\mathbb Z)$).  We claim that $y$ is not
$\Gamma(2)$-real.  Suppose otherwise.  Then $\widetilde y = \gamma y$ gives
$\widetilde \omega = \gamma'\omega$ where 
\[\gamma' = \widetilde g^{-1} \gamma g = (\widetilde g^{-1} g) g^{-1} \gamma g \in
\Gamma(2)\]
by Lemma \ref{lem-divisibility}.   But $\widetilde \omega = -\overline\omega = j\omega$ where
$j =\left(  \begin{smallmatrix} 0 & -1\\1 & 0 \end{smallmatrix} \right).$  
So $j^{-1}\gamma'$ stabilizes $\omega$ or $\gamma' = js$ for some $s\in S.$  However none of
these matrices $js$ lies in $\Gamma(2)$, which is a contradiction.  This completes the
proof of Lemma \ref{lem-upperhalfplane}.  \qed

The group $\Gamma(2)$ does not induce any identifications on the cone $iC_1$.  However $iC_1$
is preserved by the subgroup $S'= \left\{ \pm I, \pm v \right\}$ where $v = \left(
\begin{smallmatrix} 0&-1\\ 1 & 0 \end{smallmatrix}\right).$  So we obtain
\begin{thm}
The ($\tau$-)real points of $X = \Gamma(2) \backslash \mathfrak h_1$ consists of the union of
3 copies of the cone $C_1 \cong \mathbb R_{+}$,
\begin{equation*}
X_{\mathbb R} = \bigcup_{\beta \in \Gamma(2) \backslash SL_2(\mathbb Z) / S'} \beta \cdot
iC_1. \qed \end{equation*} \end{thm}
We remark that in the Baily-Borel compactification $V \cong \mathbb P^1$ these
three half-lines are glued end-to-end to form a single circle.

\section{Moduli space interpretation}\label{sec-moduli}

\subsection{}\label{subsec-PP}
Recall \cite{Lange} that a symplectic form $Q$ on $\mathbb C^n$ is {\it compatible} with the
complex structure if $Q(iu,iv) = Q(u,v)$ for all $u,v\in \mathbb C^n.$  A compatible form $Q$
is {\it positive} if the symmetric form $R(u,v) = Q(iu,v)$ is positive definite.  If $Q$ is
compatible and positive then it is the imaginary part of a unique positive definite Hermitian
form $H = R+iQ.$  Let $L\subset \mathbb C^n$ be a lattice and let $H = R+iQ$ be a positive
definite Hermitian form on $\mathbb C^n.$  A basis of $L$ is {\it symplectic} if the matrix
for $Q$ with respect to this basis is $\left( \begin{smallmatrix} 0 & I \\ -I & 0
\end{smallmatrix} \right).$  The lattice $L$ is symplectic if it admits a symplectic basis.

A principally polarized abelian variety is a pair $(A = \mathbb C^n/L, H = R+iQ)$ where $H$ is
a positive definite Hermitian form on $\mathbb C^n$ and where $L \subset \mathbb C^n$ is a
symplectic lattice relative to $Q = \text{Im}(H).$  A real structure on $(A,H)$ is a complex
anti-linear involution $\kappa:\mathbb C^n \to \mathbb C^n$ such that $\kappa(L) = L.$  A real
structure $\kappa$ on $(A,H)$ is compatible with the polarization $H$ if $Q(\kappa u, \kappa
v) = - Q(u,v)$ for all $u,v \in \mathbb C^n.$  In this case (following \cite{Silhol}) we refer
to the triple $(A, H, \kappa)$ as a {\it real principally polarized abelian variety}.  If
$(A', H', \kappa')$ is another such, then an isomorphism between them is a complex linear
mapping $\phi: \mathbb C^n \to \mathbb C^n$ such that \begin{enumerate}
\item[(a.)] $\phi(L) = L'$,
\item[(b.)] $\phi_*(H)=H'$,
\item[(c.)] $\phi_*(\kappa) = \kappa'$, that is, $\phi \kappa \phi^{-1} = \kappa'.$
\end{enumerate}

Let $(A = \mathbb C^n/L, H=R+iQ)$ be a principally polarized abelian variety.  A level $N$
structure on $A$ is a choice of basis $\left\{ U_i,V_j\right\}$ (with $ 1 \le i,j \le n$) for
the $N$-torsion points of $A,$ which is symplectic, in the sense that there exists a
symplectic basis $\left\{ u_i,v_j\right\}$ for $L$ such that
\[
U_i \equiv \frac{u_i}{N} \text{ and } V_j \equiv \frac{v_j}{N}\ \mod L\]
(for $1\le i,j \le n$).  For a given leven $N$ structure, such a choice $\left\{
u_i,v_j\right\}$ determines a mapping
\begin{equation}\label{eqn-liftmapping}
 F: \mathbb R^n \oplus \mathbb R^n \to \mathbb C^n \end{equation}
such that $F(\mathbb Z^n \oplus \mathbb Z^n) = L,$ by $F(e_i) = u_i$ and $F(f_j) = v_j$ where
$\left\{ e_i,f_j\right\}$ (with $ 1 \le i,j \le n$) is the standard basis of $\mathbb R^n
\oplus \mathbb R^n.$  The choice $\left\{u_i,v_j\right\}$ (or equivalently, the mapping $F$)
will be referred to as a {\it lift} of the level $N$ structure.  It is well defined modulo the
principal congruence subgroup $\Gamma(N)$, that is, if $F': \mathbb R^n \oplus \mathbb R^n \to
\mathbb C^n$ is another lift of the level structure, then $F' \circ F^{-1} \in \Gamma(N).$

A level $N$-structure $\left\{ U_i, V_j \right\}$ is {\it compatible} with a real structure
$\kappa$ if for some (and hence for any) lift $\left\{ u_i, v_j \right\}$ of the level
structure,
\[ \kappa({\textstyle{\frac{u_i}{N}}}) \equiv - {\textstyle{\frac{u_i}{N}}}\ \mod{L} \text{
and } \kappa({\textstyle{\frac{v_j}{N}}}) \equiv +{\textstyle{ \frac{v_j}{N}}}\ \mod L \]
for all $1 \le i,j \le n.$  In other words, the following diagram commutes $ \mod L.$
\begin{equation}\label{diag-modL} \begin{CD}
\frac{1}{N}\left( \mathbb Z^n \oplus \mathbb Z^n\right) @>>{F}> \frac{1}{N}L\\
@V{\tr{I}_-}VV @VV{\kappa}V\\
\frac{1}{N} \left(\mathbb Z^n \oplus \mathbb Z^n \right)@>{F}>> \frac{1}{N} L
\end{CD}\end{equation}
where (cf.~\S \ref{subsec-tau})
$I_- = \left( \begin{smallmatrix} -I & 0 \\ 0 & I \end{smallmatrix} \right)$

\begin{defn}\label{def-rpp}  A real principally polarized abelian variety with level $N$
structure is a quadruple $\mathcal A=\left(A = \mathbb C^n/L, H = R+iQ, \kappa,
\{U_i,V_j\} \right)$ where $(A,H,\kappa)$ is a real principally polarized abelian variety
and where $\left\{U_i,V_j\right\}$ is a level $N$-structure which is compatible with $\kappa.$
An isomorphism 
\begin{equation}\label{eqn-iso}
 \mathcal A = \left(A,H, \kappa, \{U_i,V_j\}\right) \cong 
\left(A',H',\kappa', \{ U'_i, V'_j\} \right)= \mathcal A' \end{equation}
is a complex linear mapping $\phi: \mathbb C^n \to \mathbb C^n$ such that (a.),(b.), and (c.)
of \S \ref{subsec-PP} hold, and such that \begin{enumerate}
\item[(d.)] $\phi(\frac{u_i}{N}) \equiv \frac{u'_i}{N} \ \mod{L'}$ and 
      $\phi(\frac{v_j}{N}) \equiv \frac{v'_j}{N}\ \mod{L'}$ for $1 \le i,j \le n$
\end{enumerate}
for some (and hence for any) lift $\{ u_i, v_j \}$ and $\{ u'_i, v'_j\}$ of the level
structures.
\end{defn}

\subsection{} Fix $N \ge 1.$  Each $Z \in \mathfrak {h}_n$ determines a principally polarized
abelian variety $(A_Z, H_Z)$ with level $N$ structure as follows.  Let $Q_0$ be the standard
symplectic form on $\mathbb R^n \oplus\mathbb R^n$ with matrix $\left( \begin{smallmatrix} 0 &
I \\ -I & 0 \end{smallmatrix}\right)$ (with respect to the standard basis of $\mathbb R^n
\oplus \mathbb R^n$).  Let $F_Z:\mathbb R^n \oplus \mathbb R^n \to \mathbb C^n$ be the real
linear mapping with matrix $(Z,I)$, that is,
\[ F_Z\left(\begin{smallmatrix} x \\ y \end{smallmatrix} \right ) = Zx + y.\]
Then $Q_Z=(F_Z)_*(Q_0)$ is a compatible, positive symplectic form and $L_Z = F_Z(\mathbb Z^n
\oplus\mathbb Z^n)$ is a symplectic lattice with symplectic basis $F_Z(\text{standard
basis}).$  The Hermitian form corresponding to $Q_Z$ is 
\[H_Z(u,v) = Q_Z(iu,v) + i Q_Z(u,v)= \tr{u} (\text{Im}(Z))^{-1}\bar v \text{ for } u,v \in
\mathbb C^n.\]
Then the pair $(A_Z = \mathbb C^n/L_Z, H_Z)$ is a principally polarized abelian variety.
If $z_1,z_2$ are the standard coordinates on $\mathbb C^n$ then, with respect to the above
symplectic basis of $L$, the differential forms $dz_1, dz_2$ have period matrix $(Z,I).$
If $\{e_1,\ldots,e_n, f_1,\ldots, f_n\}$ denote the standard basis elements of $\mathbb R^n
\oplus \mathbb R^n$ then the collection $\{F_Z(e_i/N), F_Z(f_i/N)\}\ \mod L$ is a level $N$
structure on $(A_Z,H_Z)$, which we refer to as the {\it standard level $N$ structure}.

Let $Z, \Omega \in \mathfrak h_n.$ Suppose $\psi: (A_{\Omega} = \mathbb C^n/L_{\Omega},
H_{\Omega}) \to (A_Z = \mathbb C^n/L_Z, H_Z)$ is an isomorphism of the corresponding
principally polarized abelian varieties, that is,  $\psi(L_{\Omega}) = L_Z$ and
$\psi_*(H_{\Omega}) = H_Z.$  Set $h = \tr{(F_Z^{-1}\psi F_{\Omega})} = \left(
\begin{smallmatrix} A& B\\ C& D \end{smallmatrix} \right).$  Then: $h\in
\mathbf{Sp}(2n,\mathbb Z)$, $\Omega = h \cdot Z$, and $\psi(M) = \tr{(CZ+D)}M$ for all $M \in
\mathbb C^n.$  Since $h\cdot Z$ is symmetric, the following diagram commutes:
\begin{equation}\label{diag-abelian1}
\begin{CD}
\left(\begin{smallmatrix} x \\ y \end{smallmatrix} \right)&\quad&
\mathbb R^n \oplus \mathbb R^n @>>{F_{h\cdot Z}}> \mathbb C^n &\quad& M \\
@VVV @VVV @VV{\psi}V @VVV \\
\tr{h}\left(\begin{smallmatrix} x \\ y \end{smallmatrix} \right) &\quad&
\mathbb R^n \oplus \mathbb R^n @>>{F_Z}> \mathbb C^n  &\quad& \tr{(CZ+D)}M
\end{CD}\end{equation}
Let $Z \in \mathfrak h_n$ with  $\widetilde{Z}= \tau(Z) = -\overline{Z}.$
Then the following diagram commutes:
\begin{equation}\label{diag-abelian2}
\begin{CD}
\left(\begin{smallmatrix} x \\ y \end{smallmatrix} \right)&\quad&
\mathbb R^n \oplus \mathbb R^n @>>{F_{Z}}> \mathbb C^n &\quad& M \\
@VVV @VVV @VVV @VVV \\
\tr{I_-}\left(\begin{smallmatrix} x \\ y \end{smallmatrix} \right) &\quad&
\mathbb R^n \oplus \mathbb R^n @>>{F_{\tau(Z)}}> \mathbb C^n  &\quad& \overline{M}
\end{CD}\end{equation}

\subsection{}
The points in $\mathfrak h_n^{\Gamma(1)}$ give rise to abelian varieties with a real
structure, as follows.  Suppose $Z \in \mathfrak h_n$, $\gamma= \left( \begin{smallmatrix} A&B
\\ C& D \end{smallmatrix}\right) \in \mathbf{Sp}(2n,\mathbb Z)$ and suppose that $\gamma \cdot
Z =\tau(Z).$  Define 
\begin{equation}
\kappa(\gamma,Z):\mathbb C^n \to \mathbb C^n\text{  by  } M \mapsto
\tr{(CZ+D)}\overline{M}.\end{equation}
Then $\kappa(\gamma,Z)$ is a real structure on $(A_Z,H_Z)$ which
is compatible with the polarization $H_Z,$ and the following diagram (which is the composition
of the two preceding diagrams) commutes:  
\begin{equation}\label{diag-kappa}
\begin{CD}
\left(\begin{smallmatrix} x \\ y \end{smallmatrix} \right)&\quad&
\mathbb R^n \oplus \mathbb R^n @>>{F_{Z}}> \mathbb C^n &\quad& M \\
@VVV @VVV @VVV @VV{\kappa(\gamma,Z)}V \\
\tr{\gamma}\tr{I_-}\left(\begin{smallmatrix} x \\ y \end{smallmatrix} \right) &\quad&
\mathbb R^n \oplus \mathbb R^n @>>{F_{Z}}> \mathbb C^n  &\quad& \tr{(CZ+D)}\overline{M}
\end{CD}\end{equation}

\begin{prop}\label{prop-gammaN}  
Let $Z \in \mathfrak h_n$ and $\gamma \in \mathbf{Sp}(2n, \mathbb R)$ and
suppose that $\widetilde{Z} = \gamma \cdot Z.$  Then $\gamma \in \Gamma(N)$ iff the real
structure $\kappa(\gamma,Z)$ on $(A_Z, H_Z)$ is compatible with the standard level $N$
structure, and in this case the quadruple
\[ \left(A_Z, H_Z, \kappa(\gamma,Z), \left\{ F_Z(e_i/N), F_Z(f_j/N) \right\} \right) \]
is a real principally polarized abelian variety with (compatible) level $N$ structure.
\end{prop}
\noindent
The proof follows immediately from the diagrams (\ref{diag-kappa}) and
(\ref{diag-modL}).  We remark that if $N \ge 3$, and if $Z \in \mathfrak h_n^{\Gamma(N)}$ then
there is a unique $\gamma \in \Gamma(N)$ such that $\widetilde{Z} = \gamma \cdot Z$ so we may
unambiguously denote $\kappa(\gamma,Z)$ simply by $\kappa_Z.$ In this case, denote by
$\mathcal A_Z$ the resulting real abelian variety with principal polarization and (compatible)
level $N$ structure.  The proof of the following theorem will appear in \S
\ref{sec-Comessatti}.

\begin{thm}\label{thm-moduli}
Fix $N = 4m$ with $m \ge 1.$  The association $Z \mapsto \mathcal A_Z$ determines a one to one
correspondence between the points of $X_{\mathbb R}$ and the set of isomorphism classes of
real principally polarized abelian varieties with level $N$ structure.
\end{thm}

\section{The Comessatti Lemma} \label{sec-Comessatti}
Recall the following theorem of Silhol \cite{Silhol} Thm. 1.4., (also \cite{Silhol2}) and
Comessatti \cite{Comessatti}.
\begin{thm}\label{thm-Silhol}  Let $(A, H, \kappa)$ be a real principally polarized abelian
variety.  Then there exists $Z = X+iY \in \mathfrak h_n$ so that $2X$ is integral, and there
exists an isomorphism of real principally polarized abelian varieties,
\[ (A, H, \kappa) \cong (A_Z, H_Z,  \sigma )\]
where $\sigma(M) = \overline{M}$ is complex conjugation.
\end{thm}
For the purposes of this paper we will need a slight restatement of this result.
\begin{lem}\label{lem-Comessatti}

{\rm (A)}  Let $Z \in \mathfrak h_n$ and suppose that $\widetilde{Z} = \gamma \cdot Z$ for
some $\gamma \in \mathbf{Sp}(2n, \mathbb Z).$  Then $Z$ is equivalent under
$\mathbf{Sp}(2n,\mathbb Z)$ to an element $X+iY \in \mathfrak h_n$ such that $2X$ is integral. 

{\rm (B)}  Let $(A,H,\kappa)$ be a real principally polarized abelian variety.    Then there
exists $\gamma \in \mathbf{Sp}(2n, \mathbb Z)$ and there exists $Z = X + iY \in \mathfrak h_n$
such that $\widetilde{Z} = \gamma \cdot Z$ and $2X$ is integral; and there exists an
isomorphism of real principally polarized abelian varieties
\[ (A,H,\kappa) \cong (A_Z, H_Z, \kappa(\gamma,Z)).\]
\end{lem}

\subsection{Proof}
For part (B), Theorem \ref{thm-Silhol} provides $Z = X+iY \in \mathfrak h_n$ so that $2X$ is
integral.  Take $\gamma = \left( \begin{smallmatrix} I & -2X \\ 0 & I \end{smallmatrix}
\right).$  Then $\widetilde{Z} = - \overline{Z} = \gamma \cdot Z$ and $\kappa(\gamma,Z)(M) =
\overline{M}=\sigma(M)$ by (\ref{diag-kappa}).  

For part (A), given $Z$ and $\gamma$, apply part (B) to
the real abelian variety $(A_Z, H_Z, \kappa(\gamma,Z))$ to obtain an isomorphic real abelian
variety $(A_{Z'}, H_{Z'}, \kappa(\gamma',Z'))$ such that $Z' = X'+iY'$ with $2X'$ integral.
Then $(A_Z,H_Z)$ and $(A_{Z'},H_{Z'})$ are isomorphic principally polarized abelian varieties
so there exists $g \in \mathbf{Sp}(2n,\mathbb Z)$ with $Z' = g\cdot Z.$  \qed

\subsection{Proof of Proposition \ref{add}}
By Comessatti's lemma, $Z$ is equivalent (via some  $h\in
\mathbf{Sp}(2n,\mathbb Z)$) to some element $X+iY \in \mathfrak h_n$ with $2X \in M_{n\times
n}(\mathbb Z).$  We claim, in this case, that $X$ is integral.
Translation by $X$ is given by the symplectic matrix 
\begin{equation*}T_X=\bigl(
\begin{smallmatrix} I & X \\ 0 & I \end{smallmatrix} \bigr)
\end{equation*}
so we may write $Z = h(X+iY)=hT_XiY.$  Then $\tau Z = \tau(hT_X)iY = \gamma Z,$ so the
following element $\gamma^{-1}\tau(hT_{X})(hT_X)^{-1}$
fixes $Z$.  By our assumption on $Z$, this implies that 
\begin{equation}\label{eqn-1}
\gamma = \pm I \tau(hT_X) (hT_X)^{-1}
\end{equation} or, 
\begin{equation*}
\pm T_{-2X} = \tilde h^{-1}\gamma h = (\tilde h^{-1}h) (h^{-1}\gamma h) \in
\Gamma(2).\Gamma(4m) 
\end{equation*}
(using Lemma \ref{lem-divisibility} and the fact that $\Gamma(4m)$ is normal in
$\mathbf{Sp}(2n,\mathbb Z)$).  So $2X$ is ``even'', hence $X$ is integral.

If the plus sign occurs in (\ref{eqn-1}) then $\gamma = \tau(hT_X)(hT_X)^{-1}$ hence
$\mathfrak h_n^{\gamma} = hT_X \cdot iC_n.$  If the minus sign occurs, set $\omega =
\left(\begin{smallmatrix} 0 & -I\\ I & 0 \end{smallmatrix} \right).$  Then $\gamma =
\tau(hT_X\omega)(hT_X\omega)^{-1}$ and hence $\mathfrak h_n^{\gamma} = hT_X\omega\cdot iC_n.$
\qed

In the next proposition we strengthen these results to include level structures.

\begin{prop}\label{prop-thegamma}
Suppose $N \ge 3.$  Suppose $\mathcal A$ is a real abelian variety with principal polarization
and level $N$ structure.  Then there exists $Z \in
\mathfrak h_n$, there exists $\gamma \in \Gamma(N)$ such that $\gamma \cdot Z = \tau(Z)$, and
there exists an isomorphism 
\[
\phi: \mathcal A_Z \to \mathcal A
\]
of real principally polarized abelian varieties with level structures.  If $N=4m$ (and $m \ge
1$) the cohomology class $[f_{\gamma}] \in H^1(\mathbb C/
\mathbb R, \Gamma(4m))$ is uniquely determined by the isomorphism class of $\mathcal
A.$\end{prop}

\subsection{Proof}  Write $\mathcal A = (A= \mathbb C^n/L, H= R + iQ, \kappa, \{
U_i,V_j\}).$  By Lemma \ref{lem-Comessatti}, there exists $Z' \in
\mathfrak h_n$, there exists $\gamma' \in \mathbf{Sp}(2n,\mathbb Z)$ with $\gamma'\cdot Z' =
\tau(Z')$, and there exists $\phi':\mathbb C^n \to \mathbb C^n$ such that $\phi'$ induces an
isomorphism  $(A,H,\kappa)\cong (A_{Z'}, H_{Z'}, \kappa(\gamma',Z')),$  that is:
\begin{enumerate}
\item[(a.)] $\phi'(L_{Z'}) = L$
\item[(b.)] $\phi'_*(H_{Z'}) = H$
\item[(c.)] $\phi'_*(\kappa(\gamma',Z')) = \kappa$.
\end{enumerate}
The isomorphism $\phi'$ must be modified because it does not necessarily take the
standard level $N$ structure on $(A_{Z'}, H_{Z'})$ to the given level $N$ structure on
$(A,H).$  Choose a lift $\left\{ u_i, v_j \right\}$ ($1 \le i,j \le n$) of the level $N$
structure on $(A,H)$ and let $F: \mathbb R^2 \oplus \mathbb R^2 \to \mathbb C^2$ be the
resulting mapping (\ref{eqn-liftmapping}).  Define 
\begin{align}
 \tr{g}^{-1} &= F^{-1} \circ \phi' \circ F_{Z'} \in \mathbf{Sp}(2n, \mathbb Z) \\
Z &= g \cdot Z'  \\
\label{eqn-gammaprime}
\gamma &= \tau(g) \gamma' g^{-1} = I_- g I_- \gamma' g^{-1}.\end{align}  If $g = \left(
\begin{smallmatrix} A&B\\C&D \end{smallmatrix}\right)$ define $\psi:\mathbb C^n \to \mathbb
C^n$ by $\psi(M) = \tr{(CZ+D)}M.$  We claim that $\gamma \in \Gamma(N),$ that $\widetilde{Z} =
\gamma \cdot Z$ and that the mapping
\[\phi = \phi' \circ \psi:\mathbb C^n \to
\mathbb C^n\]
induces an isomorphism $\mathcal A_Z \cong \mathcal A$ between principally polarized real
abelian varieties with (compatible) level $N$ structures.

The matrix $g$ is defined so that the bottom square in the following diagram commutes; 
by (\ref{diag-abelian1}) the top square also commutes.  The mapping $\phi$ is the composition
down the right-hand column.
\begin{equation} \label{diag-abelian1prime}
\begin{CD}
\mathbb R^n \oplus \mathbb R^n @>{F_{Z}}>> \mathbb C^n  \\
@V{\tr{g}}VV @VV{\psi}V  \\
\mathbb R^n \oplus \mathbb R^n @>>{F_{Z'}}> \mathbb C^n \\
@V{\tr{g}^{-1}}VV @VV{\phi'}V \\
\mathbb R^n \oplus \mathbb R^n @>>{F}> \mathbb C^n
\end{CD}\end{equation}
Clearly, $\widetilde Z = \gamma \cdot Z,$ $\phi_*(L_Z) = L$ and $\phi_*(H_Z) = H.$  Now let us
check that
\begin{equation}\label{eqn-phikappa}
\phi_*(\kappa(\gamma,Z)) = \kappa.\end{equation}  
By \S 7.6 (c.) it suffices to check that $\psi \kappa(\gamma,Z) \psi^{-1} =
\kappa(\gamma',Z').$  But this follows from direct calculation using
$\kappa(\gamma,Z) = F_Z \tr{\gamma} \tr{I}_{-} F_Z^{-1},$ $\kappa(\gamma',Z') = F_{Z'}
\tr{\gamma'} \tr{I}_{-} F_{Z'}^{-1},$ $\psi = F_{Z'} \tr{g} F_Z^{-1},$ and
(\ref{eqn-gammaprime}).

By (\ref{diag-abelian1prime}), $F = \phi\circ F_Z.$ Therefore $\phi$ preserves the level
structures. Since the given level structure on $(A,H)$ is compatible with $\kappa$, equation
(\ref{eqn-phikappa}) guarantees that the standard level $N$ structure on $(A_Z,H_Z)$ is
compatible with $\kappa(\gamma,Z).$  It follows from Proposition \ref{prop-gammaN} that
$\gamma \in \Gamma(N).$    In summary, we
have shown that $(A_Z,H_Z, \kappa(\gamma,Z), \left\{F_Z(e_i/N), F_Z(f_j/N)\right\})$ is a real
principally polarized abelian variety with (compatible) level $N$ structure, and that the
isomorphism $\phi$ preserves both the real structure and the level structure.  

Finally, an isomorphism $\phi: \mathcal A \to \mathcal A'$ (\ref{eqn-iso}) between real
principally polarized abelian varieties with level
$N = 4m $ structures determines a cohomological equivalence between the corresponding
1-cocycles as follows.  Choose $Z, Z' \in \mathfrak h_n$ and $\gamma, \gamma' \in \Gamma(4m)$ 
for $\mathcal A$ and $\mathcal A'$ respectively, as above, and set
\[ g = {}^{t}(F_Z^{-1} \circ \phi \circ F_{Z'}).\]
Then $ g \in \Gamma(4m)$ by \S \ref{def-rpp}(d.).  Since $\phi \kappa \phi^{-1} = \kappa'$,
diagrams \ref{diag-abelian1} and \ref{diag-kappa} give
\[ \tr{g} \tr{\gamma} \tr{I_-} \tr{g^{-1}} = \tr{\gamma'}I_{-} \]
or $\gamma' = \tilde g^{-1} \gamma g.$  Hence $[f_{\gamma}] = [f_{\gamma'}]$ in $H^1(\mathbb
C/ \mathbb R, \Gamma(4m)).$
\qed

\subsection{Proof of Theorem \ref{thm-moduli}}
Let $\mathcal I$ denote the moduli space of isomorphism classes of real abelian varieties with
principal polarization and level $N$ structure.  Let
\[ S = \bigcup_{\gamma \in \Gamma(N)} \mathfrak h_n^{\gamma} \subset \mathfrak h_n\]
denote the set of $\Gamma(N)$-real points in $\mathfrak h_n$ and let $\Gamma(N) \backslash S$
denote its image in $X = \Gamma(N) \backslash \mathfrak h_n.$ Let $\Phi: S \to \mathcal I$ be
the association 
\[Z \mapsto (A_Z, H_Z, \kappa_Z, \{F_Z(e_i/N), F_Z(f_j/N)\}).\] Here, $\kappa_Z:
\mathbb C^n \to \mathbb C^n$ is the real structure given by (\ref{diag-kappa}), that is, 
$\kappa(M) = \tr{(CZ+D)}\overline{M}$ where $\gamma = \left( \begin{smallmatrix} A&B\\C&D
\end{smallmatrix}\right)\in \Gamma(N)$ is the unique element such that $\gamma \cdot Z = -
\overline{Z}.$  Then Proposition  \ref{prop-thegamma} says that $\Phi$ is surjective.  If $Z,
Z' \in S$ and if $Z' = \gamma \cdot Z$ for some $\gamma \in \Gamma(N)$ then by
(\ref{diag-abelian1}) the real abelian variety $(A_Z, H_Z, \kappa_Z, \{F_Z(e_i/N),
F_Z(f_i/N)\})$ is isomorphic to $(A_{Z'}, H_{Z'}, \kappa_{Z'}, \{F_{Z'}(e_i/N),
F_{Z'}(f_i/N)\}).$  Therefore $\Phi$ passes to a mapping $\Gamma(N) \backslash S \to \mathcal
I.$  

On the other hand, if $Z, \Omega \in S$ and if there is an isomorphism $\psi: \mathcal
A_{\Omega} \to \mathcal A_Z$ then by (\ref{diag-abelian1}) there exists $h \in
\mathbf{Sp}(2n,\mathbb Z)$ such that $\Omega = h \cdot Z.$  Since the mapping $\psi$ also
preserves the level structures, it follows also from (\ref{diag-abelian1}) that $h \in
\Gamma(N).$  Hence the mapping $\Gamma(N) \backslash S \to \mathcal I$ is
also injective.  By Theorem \ref{thm-real}, the quotient $\Gamma(N)\backslash S$ is precisely
the variety $X_{\mathbb R}$ of real points in $X$.  \qed

\section{Baily-Borel Compactification}\label{sec-BB}
  Let $\overline{\mathfrak h}_n$ denote the Satake partial compactification of $\mathfrak h_n$
which is obtained by attaching all rational boundary components, with the Satake topology.
The group $Sp(2n,\mathbb Q)$ acts on $\overline{\mathfrak h}_n$, the involution
$\tau:\mathfrak h_n\to \mathfrak h_n$ extends to $\overline{\mathfrak h}_n$, and $\tau(gx) =
\tau(g)\tau(x)$ for any $g\in Sp(2n,\mathbb Q)$ and every $x\in \overline{\mathfrak h}_n.$
Denote by $\pi:\overline{\mathfrak h}_n \to \Gamma(4m)\backslash \overline{\mathfrak h}_n =
V$ the projection to the Baily-Borel compactification of $X$.  The involution
$\tau$ passes to complex conjugation $\tau:V \to V$, whose fixed points we denote by
$V_{\mathbb R}.$  Clearly $X_{\mathbb R} \subset V_{\mathbb R}.$  Define
$\overline{X}_{\mathbb R}$ to be the closure of $X_{\mathbb R}$ in $V_{\mathbb R}.$

In \cite{Shimura}, Shimura shows
that the $\Gamma(N)$-automorphic forms on $\mathfrak h_n$ are generated by those automorphic
forms with rational Fourier coefficients.  It follows that: 
\begin{thm}  There exists a natural rational structure on the Baily-Borel compactification
$V$ of $X$ which is compatible with the $\tau$-real structure.
\end{thm}

The Baily-Borel compactification $V$ is stratified by finitely many strata of the form
$\pi(F)$ where $F \subset \overline{\mathfrak h}_n$ is a rational boundary component.  Each
such $F$ is isomorphic to some Siegel space $\mathfrak h_k$, in which case we say the stratum
$\pi(F)$ is a boundary stratum of rank $k$.  Let $V^r$ denote the union of all boundary strata
of rank $n-r.$  In Proposition \ref{prop-BBcompact} we will prove that
\begin{equation}\label{eqn-BBcompact}
V^0_{\mathbb R} \cup V^1_{\mathbb R} \subset \overline{X}_{\mathbb R} \subset V_{\mathbb R}
\end{equation}
where $V^r_{\mathbb R} = V^r \cap V_{\mathbb R}.$

As in \S \ref{sec-Galois}, if $F$ is a rational boundary component of $\mathfrak h_n$
which is preserved by the involution $\tau$, and if $\gamma \in \Gamma(4m)$, we define the set
of $\gamma$-real points of $F$ to be
\[ F^{\gamma} = \left\{ x \in F |\ \tau(x) = \gamma x \right\}. \]
Then $\pi(F^{\gamma}) \subset V_{\mathbb R}.$  A {\it $\Gamma(4m)$-real boundary pair}
$(F,\gamma)$ (of rank $q$) consists of a rational boundary component $F$ (of rank $q$) and an
element $\gamma \in \Gamma(4m)$ such that $F^{\gamma} \ne \phi.$ (Hence $\widetilde{F} =
\gamma F.$)  We say two real boundary pairs $(F,\gamma)$ and $(F_1,\gamma_1)$ are
equivalent if the resulting locus of real points $\pi(F^{\gamma}) = \pi(F_1^{\gamma_1})$
coincide.  If $(F,\gamma)$ is a real boundary pair and if $g\in \Gamma(4m)$ then
$(gF,\tilde g \gamma g^{-1})$ is an equivalent real boundary pair. 

\subsection{Standard boundary components}\label{subsec-standardboundary}
Fix an integer $q$ with $ 1 \le q < n.$  The Siegel upper halfspace
\begin{equation*}
\mathfrak h_{q} = \{ Z\in M_{q\times q}(\mathbb C) | \ \tr Z=Z 
\text{ and } \text{Im}(Z)>0\} \end{equation*}
is attached to $\mathfrak h_n$ as a limit of matrices in $M_{n\times n}(\mathbb C)$ by 
\[Z \mapsto \lim_{Y} \left(\begin{smallmatrix} Z&0\\0&iY
\end{smallmatrix}\right).\]
See (\ref{eqn-iCn}).  Here, $Y\in C_{n-q}$ is a positive definite symmetric matrix of order
$n-q$, and the limit is taken as $Y \to \infty C_{n-q}.$  This means that for any $Y_0 \in
C_{n-q}$, the difference $Y-Y_0$ is eventually contained in $C_{n-q}$ (or, alternatively, that
all the eigenvalues of $Y$ converge to $\infty$).    
Denote this mapping by $\phi: \mathfrak h_{q} \to \overline{\mathfrak h}_n.$   Its image
$F_{q} = \phi(\mathfrak h_{q}) \subset \overline{\mathfrak h}_n$ is called the {\it
standard boundary component} of (maximal) rank $q.$   The normalizer in $Sp(2n,\mathbb Z)$
of $F_{q}$ is the parabolic subgroup   
\begin{equation*}
P_{q} =  \left( \begin{array}{cc|cc}
 A& 0& B& *\\
 *&*& *& *\\ \hline
 C&0& D& * \\
 0& 0& 0&* \end{array} \right)
\end{equation*}
with unipotent radical
\begin{equation*}
\mathcal U(P_q) = \left( \begin{array}{cc|cc}
I_q & 0 & 0 & b \\
a &I_{n-q} & \tr b & d \\ \hline
0&0&I_q & -\tr a \\
0&0&0&I_{n-q}
\end{array}\right) \end{equation*}
and Levi factor $L(P_{q}) = G_hG_{\ell}$ with
\begin{equation}\label{eqn-splitting}
G_h = \left( \begin{array}{cc|cc}
A&0&B&0 \\ 0&I_{n-q} &0&0 \\ \hline
C&0&D&0 \\ 0&0&0&I_{n-q}
\end{array} \right)\text{ and } 
G_{\ell} = \left( \begin{array}{cc|cc}
I_q &0&0&0 \\ 0 & T & 0&0\\ \hline
0&0&I_q &0 \\ 0&0&0& {}^tT^{-1}
\end{array} \right) \end{equation}
being the ``Hermitian'' and ``linear'' factors respectively, where $\left( \begin{smallmatrix}
A&B\\C&D \end{smallmatrix}\right) \in
Sp(2m,\mathbb R)$ and $T\in \text{GL}_{n-q}(\mathbb R).$  The subgroup $\mathcal
U(P_{q})G_{\ell}$ is normal in $P_{q}$ and we denote by 
\[ \nu:P_{q} \to G_h\cong \mathbf{Sp}(2q, \mathbb R)\] the projection to
the quotient.  Then $\nu$ commutes with the involution $\tau.$

The boundary component $F_{q}$ is preserved by $\tau$.  The set 
\begin{equation*}
F_{q}^{\tau}=\phi(iC_{q}) = \left\{ \phi(iY)|\  Y >0\right\}
\end{equation*}
of $I$-real points in $F_{q}$ is just the set of $\tau$-fixed points in $F_{q}$ and it may
be canonically identified with the cone of positive definite matrices of order $q.$  Denote
by $iI_{q}$ its canonical basepoint.  The boundary component $F_{q}$ is attached to
$\mathfrak h_n$ so that this cone $\phi(iC_{q})$ is contained in the closure of the cone
$iC_n.$

\begin{prop}\label{lem-noB}
  Let $(F,\gamma)$ be a $\Gamma(4m)$-real boundary pair of rank $q.$  Then there exists $a
\in \mathbf{Sp}(2n,\mathbb Z)$ so that $a(F_q) = F$ and
\begin{equation}\label{eqn-firsttry}
\tilde a^{-1} \gamma a = \left( \begin{matrix} A & B \\ 0 & \trA^{-1} \end{matrix}
\right) \in \ker(\nu). \end{equation}
Moreover, we may take $B=0$, that is, there exists $\gamma' \in \Gamma(4m)$ and $g \in
\mathbf{Sp}(2n,\mathbb Z)$ so that $F^{\gamma'} = F^{\gamma}$, $g(F_q) = F$,
and so that \begin{equation}\label{eqn-secondtry}
 \tilde g^{-1} \gamma' g = \left( \begin{matrix} A & 0 \\ 0 & \trA^{-1}\end{matrix}
\right)\in \ker(\nu).\end{equation}
\end{prop}

\subsection{Proof}\label{pf-noB}  Choose $b \in \mathbf{Sp}(2n,\mathbb Z)$ so that $b F_q =
F.$  Then $w = \tilde b^{-1} \gamma b = \tilde b^{-1} b \cdot b^{-1} \gamma b \in
\Gamma(2).$  Since it also preserves the standard boundary component $F_q,$ we have:  $w
\in P_q$.  

The Hermitian part $\nu(w)\in \mathbf{Sp}(2n, \mathbb Z)$ acts on the standard boundary
component $F_q \cong \mathfrak h_q.$  By Lemma \ref{lem-divisibility},
$\nu(\widetilde{w})\nu(w) \in
\Gamma(4)$ which is torsion-free.  By assumption, the set of $\nu(w)$-real points
$F_q^{\nu(w)}$ is nonempty.  Choose a point $x \in F_q^{\nu(w)}$ whose stabilizer
in $\mathbf{Sp}(2m,\mathbb Z)$ consists of $\pm I:$  by Proposition \ref{lem-Fgamma} such
points exist and are even dense in $F_q^{\nu(w)}.$    Then, by Proposition \ref{add} there
exists $h \in \mathbf{Sp}(2m,\mathbb Z)$ such
that $\tilde h h^{-1} = \nu(w)$ and $F_q^{\nu(w)} = h(iC_q).$  Let us identify the
element $h$ with its image in $\mathbf{Sp}(2n, \mathbb Z)$ as in (\ref{eqn-splitting}) and set 
$v = \tilde h^{-1} w h.$  Then $\nu(v) = 1.$  The following diagram may help to sort out
these various transformations.
\begin{equation*}
\begin{CD}
F_q @>{\tilde h}>> F @>{\tilde b}>> F\\
@A{v}AA @AA{w}A @AA{\gamma}A \\
F_q @>>{h}> F @>>{b}> F \end{CD}\end{equation*}
Then the element $a = bh \in \mathbf{Sp}(2n,\mathbb Z)$ has the desired properties, that is,
$a\cdot F_q = F$ and $\nu(v)=1$ where  $v = \tilde a^{-1} \gamma a \in \Gamma(2).$

Now let us prove the ``moreover'' part of Proposition \ref{lem-noB}.
  By Lemma \ref{lem-divisibility}, both
$\tilde \gamma \gamma^{-1}$ and $\tilde \gamma \gamma$ are in $\Gamma(8m).$  Then $\tilde
\gamma  \gamma = \tilde a \tilde v  v \tilde a^{-1} \in \Gamma(8m),$ hence $\tilde v  v \in
\Gamma(8m).$ Calculating  $\tilde v v \equiv I \mod{8m}$ gives 
\begin{equation}\label{eqn-ABmod8}
AB \equiv B \trA^{-1}  \mod{8m} \text{ and } A^2 \equiv I \mod{8m}.\end{equation}  
Since $v =\tilde a^{-1}a a^{-1}\gamma a \in
\Gamma(2)$ the matrix $A$ is integral and $B$ is even, so 
\[ x = \left( \begin{matrix} I & -\frac{1}{2}A^{-1}B \\ 0 & I \end{matrix} \right) \]
is integral.  Since $v \in \ker(\nu)$ we see that $A =
\left(\begin{smallmatrix}*&0\\ *&*\end{smallmatrix}\right)$ and $B =
\left(\begin{smallmatrix}0&*\\ *&*\end{smallmatrix}\right)$ from which it follows that  $x \in
\mathcal U(P_q) \subset \ker(\nu).$  Set $g = ax$
and set $u' = \tilde x^{-1} v x.$  Direct computation with the matrices for $x$ and $v$ gives 
\[u' = \left( \begin{matrix} A & B' \\ 0 & \trA^{-1}\end{matrix} \right)\]
where $B' = \frac{1}{2}B - \frac{1}{2}A^{-1}B\trA^{-1}.$  Using (\ref{eqn-ABmod8}) gives $B'
\equiv 0 \mod{4m}.$  The following diagram may help to explain these transformations.
\begin{equation*}\begin{CD}
F_q @>{\tilde x}>> F_q @>{\tilde a}>> F\\
@A{u'}AA @AA{v}A @AA{\gamma}A \\
F_q @>>{x}> F_q @>>{a}> F
\end{CD}\end{equation*}
Now decompose $u' = uu_2$ where
\[ u = \left( \begin{matrix} A & 0 \\ 0 & \trA^{-1}
\end{matrix} \right)\text{ and  }u_2 = \left( \begin{matrix} I & A^{-1}B' \\ 0 & I
\end{matrix} \right) \in \ker(\nu).\]  
Set $\gamma' = \tilde g u g^{-1}.$  We will verify that $g$ and $\gamma'$ satisfy the
conclusions of Lemma \ref{lem-noB}.

First note that $u_2 \in \Gamma(4m)$ since $B' \equiv 0 \mod{4m}.$  Then $\gamma' \in
\Gamma(4m)$ because
\[\gamma = \tilde a v a^{-1} =  \tilde g u' g^{-1} = (\tilde g u g^{-1})(g u_2 g^{-1})
 = \gamma'(gu_2g^{-1})\in\Gamma(4m).\]
We have already verified that $gF_q =F$ and that $u=\tilde g^{-1} \gamma' g$ has the desired
form.  Since $u_2$ acts trivially on $F_q$ we see that $y \in F^{\gamma'}$ iff
\[ \tilde y = \gamma ' y = \tilde g u g^{-1} y = \tilde g u u_2 g^{-1} y = \gamma y\]
iff $y \in F^{\gamma}.$  Hence $F^{\gamma'} = F^{\gamma}.$  \qed

\subsection{Some nearby boundary components}\label{sec-opposite}
Fix $r$ with $ q\le r \le n$ and set $s = n-r.$  Define 
\begin{equation*} j_r = \left(
\begin{array}{cc|cc}
I_{r} & 0&0 &0 \\ 0&0 &0 & I_s \\ \hline
0&0 & I_{r} & 0\\0 & -I_s &0 & 0\end{array} \right)
\text{ so } j_r^2 = \left(
\begin{array}{cc|cc}
I_r &0&0&0\\0&-I_s&0&0\\ \hline 0&0&I_r&0\\0&0&0&-I_s
\end{array}\right).
\end{equation*}
Then $\tilde j_r = j_r^{-1}$.  Although $j_r \notin P_{q},$ its square $j_r^2$ preserves
$F_{q}$ and in fact it acts as the identity on $F_{q}.$  Let $E_{q,r} = j_r(F_{q}) =
j_r^{-1}(F_{q})$ and let $Q_{q,r}$ be its normalizing maximal parabolic subgroup:  it is
$j_r$-conjugate to $P_{q}.$  (If $r=n$ then $j_r$ is the identity and $E_{q,r} = F_q$).  The
involution $\tau$ preserves the boundary component $E_{q,r}$ and {\it the set of $\tau$-fixed
points
\[E_{q,r}^{\tau}= \left\{ x \in E_{q,r} \left| \right. \tilde x = x \right\} = j_r\cdot
F_{q}^{\tau} \]
is contained in the closure of $\mathfrak h_n^{\tau}=iC_n$} as follows.   Although $j_r$ does
not preserve the cone $iC_n$, it does preserve the sub-cone $S_{q,r}\subset iC_n$ of elements
\[ \left( \begin{matrix} iY_1 & 0&0 \\ 0 & iY_2&0 \\ 0&0&iY_3 \end{matrix} \right)\in
\mathfrak h_n \]
where $Y_1\in C_{q},$  $Y_2\in C_{r-q},$ and $Y_3 \in C_s.$    Moreover, $F_{q}^{\tau}$
is contained in the closure of $S_{q,r}$.  Therefore \begin{equation}\label{eqn-EF}  
 E_{q,r}^{\tau} = j_rF_{q}^{\tau} \subset j_r\overline{S_{q,r}} = \overline{S_{q,r}} \subset
\overline{iC_n}.\end{equation}

\begin{prop}\label{prop-killer}
Let $(F, \gamma)$ be a $\Gamma(4m)$-real boundary pair of rank $q$.  Let $g \in
\mathbf{Sp}(2n,\mathbb Z).$  Suppose that $g(F_q) = F$ and that $u = \tilde g^{-1} \gamma g
\in \ker(\nu).$  Suppose also that there exists $r$ with $q \le r \le n$ so that 
\[ \tilde j_r^{-1} u j_r = j_r u j_r \in \Gamma(4m).\]
Define $\omega = \tau(g j_r) (gj_r)^{-1}.$  Then $\omega \in \Gamma(4m)$ and
\[ F^{\gamma} = F^{\omega} \subset \overline{\mathfrak h_n^{\omega}}. \]
Consequently the resulting set $\pi(F^{\gamma})$ of real points is contained in the closure
$\overline{X}_{\mathbb R}.$
\end{prop}

\subsection{Proof}
Calculate 
\begin{equation*}
\omega=\tau(gj_r) (gj_r)^{-1} = \gamma \left((gj_r)(j_ruj_r)^{-1}(gj_r)^{-1}\right) \in
\Gamma(4m)
\end{equation*}
(using the fact that $\tilde j_r = j_r^{-1}$) which proves the first statement.  Since $u \in
P_q$ acts trivially on $F_q$, the same is true of $j_r^2u.$  Hence,  $x \in F^{\omega}$ if and
only if 
\[ \tilde x = \omega x = \tau(gj_r)j_r^{-1}(g^{-1}x) = \tau(gj_r)j_r^{-1}j_r^2u(g^{-1}x)
= \tau(gj_r)\tau(gj_r)^{-1}\gamma x = \gamma x\]
which holds if and only if $x \in F^{\gamma}.$  The following diagram may help in placing
these elements,
\begin{equation*} \begin{CD}
E_{q,r} @>{\tilde j_r}>> F_q @>{\tilde g}>> F \\
@AAA @A{u}AA @AA{\gamma}A \\
E_{q,r} @>>{j_r}> F_q @>>{g}> F
\end{CD}\end{equation*}
We claim  that $F^{\omega} = gj_r(E_{q,r}^{\tau}).$  In fact, $x \in E_{q,r}^{\tau}$ if and
only if $\tilde x = x$ which holds iff
\begin{equation*}
\tau(gj_r)\tilde x = \tau(gj_r)x=\tau(gj_r)(gj_r)^{-1}(gj_r)x = \omega (gj_r)x \end{equation*}
which holds iff $ (gj_r)x \in F^{\omega}.$  Similarly, $ gj_r(\mathfrak h_n^{\tau}) =
\mathfrak h_n^{\omega}.$  Hence
\[ F^{\omega} = gj_r(E_{q,r}^{\tau}) \subset gj_r\overline{\mathfrak h_n^{\tau}} =
\overline{\mathfrak h_n^{\omega}}  \]
using (\ref{eqn-EF}).  \qed

\section{Corank one strata}
\subsection{}
In this section we specialize to the case $q=n-1$, that is, we consider only
the boundary strata of maximal rank.  As in \S \ref{sec-BB}, $F_{n-1}$ denotes the standard
boundary component of corank 1, $P_{n-1}$ denotes its normalizing parabolic subgroup and
$G_h$ and $G_{\ell}$ refer to the Hermitian and
linear factors (\ref{eqn-splitting}) of the Levi quotient $L(P_{n-1}).$

\begin{thm}\label{prop-BBcompact}
Let $F$ be a proper rational boundary component of $\mathfrak h_n$ with (maximal) rank $n-1$,
let $\gamma \in \Gamma(4m)$, and suppose that $F^{\gamma} \ne \phi.$  Then $F^{\gamma}$ is
contained in the closure of the set $\mathfrak h_n^{\Gamma(4m)}$ of $\Gamma(4m)$-real points
of $\mathfrak h_n.$
\end{thm}

\subsection{Proof} By Proposition \ref{lem-noB}, there exists $\gamma' \in \Gamma(4m)$ and $g
\in \mathbf{Sp}(2n,\mathbb Z)$ so that
$gF_{n-1} = F$, $F^{\gamma} = F^{\gamma'}$ and so that 
\[ u = \tilde g^{-1} \gamma' g = \left( \begin{matrix} A & 0 \\ 0 & \trA^{-1} \end{matrix}
\right)\in \ker(\nu).  \]
This implies that $A = \left( \begin{smallmatrix} I_{n-1} & 0 \\ a & \pm 1 \end{smallmatrix}
\right).$  Note that $\tilde u u = g^{-1} \tilde \gamma \gamma g \in \Gamma(8m)$ by Lemma
\ref{lem-divisibility}.  Hence $a \equiv 0 \mod{4m}.$

If the plus sign occurs then this says that $u \in \Gamma(4m).$  Let $\omega = \tilde g
g^{-1}.$  Then Proposition \ref{prop-killer} (with $q = n-1$ and $r = n$) implies that
$F^{\gamma'} \subset \overline{\mathfrak h_n^{\omega}}.$

If the minus sign occurs then 
\begin{equation*} j_{n-1}u j_{n-1} = \left( \begin{array}{cc|cc}   
I&0 &0 & 0 \\0 & 1 & 0 & 0 \\ \hline 
0 & -\tr a &I & 0 \\ -a & 0 & 0 & 1
\end{array} \right) \in \Gamma(4m)\end{equation*}
so we may apply Proposition \ref{prop-killer} (with $q=r = n-1$) to conclude that $F^{\gamma'}
\subset \overline{\mathfrak h_n^{\omega}}$ where $\omega =  \tau(g j_{n-1}) (gj_{n-1})^{-1}.$
\qed

\section{The principal congruence group $\Gamma({2^k})$}\label{sec-2k}
\subsection{}  Throughout this section we let $\Gamma = \Gamma({2^k}) \subset \mathbf{Sp}(2n,
\mathbb Z)$ be the principal congruence subgroup of level ${2^k}$ with $k \ge 2.$  As in \S
\ref{sec-BB}, let $X = X_{\mathbb C}= \Gamma \backslash \mathfrak h_n$ and let $V =V_{\mathbb
C}= \Gamma \backslash \overline{\mathfrak h}_n$ be its Baily-Borel compactification with
projection $\pi: \overline{\mathfrak h}_n \to V.$  Let $\overline{X}_{\mathbb R}$
be the closure of $X_{\mathbb R}$ in $V.$  In this section we will prove that
$\overline{X}_{\mathbb R} = V_{\mathbb R}.$  The proof of the following lemma will appear in
\S \ref{subsec-proofoflemma}.

\begin{lem}\label{lem-GLconj}  Let $A \in \mathbf{GL}(n, \mathbb Z).$  Suppose that 
\begin{equation*}
A \equiv I \mod 2 \text{ and } A^2 \equiv I \mod {2^{k+1}}.\end{equation*}
Then there exists $p \in \mathbf{GL}(n,\mathbb Z)$ so that
\begin{equation}\label{eqn-pm1}
 p^{-1}Ap \equiv \left( \begin{smallmatrix} \pm 1 &&&&\\ & \pm 1 &&& \\ && \cdots && \\
&&& \pm 1 & \end{smallmatrix} \right) \mod {2^k}. \end{equation}
Moreover if the matrix of $A$ with respect to the standard basis of $\mathbb Z^n$ is $\left(
\begin{smallmatrix} I_q & 0 \\ * & * \end{smallmatrix}\right)$ then it is possible to choose
$p$ to be of the form $p = \left( \begin{smallmatrix} I_q & 0 \\ * & * \end{smallmatrix}
\right).$
\end{lem}

\begin{lem}\label{lem-IrIs}
Let $(F,\gamma)$ be a real boundary pair of rank $q.$  Then there exists $r \ge q$,
there exists $\gamma' \in \Gamma({2^k})$ and there exists  $g\in \mathbf{Sp}(2n,\mathbb
Z)$ such that$F^{\gamma} = F^{\gamma'}$,  $gF_q = F$,  $\tilde g^{-1} \gamma' g \in
\ker(\nu),$ and so that
\begin{equation}\label{eqn-normalform} \tilde g ^{-1} \gamma' g \equiv \left(
\begin{array}{cc|cc}
I_r &0 & 0 & 0 \\
0 & -I_s &0 &0 \\ \hline
 0 & 0 & I_r & 0 \\ 0 & 0 & 0 & -I_s \end{array} \right) \mod {2^k} \end{equation}
where $s = n-r.$
\end{lem}

\subsection{Proof} By Proposition \ref{lem-noB} there exists $\gamma' \in
\Gamma({2^k})$ and there exists $a \in \mathbf{Sp}(2n,\mathbb Z)$ so that $F^{\gamma} =
F^{\gamma'}$, $aF_q = F$, and 
\[u=\tilde a^{-1}\gamma' a = \left( \begin{matrix} A & 0
\\ 0 & \trA^{-1} \end{matrix} \right)\in \ker(\nu).\]  
Then $\tilde \gamma' \gamma' = a u^2 a^{-1} \in \Gamma({2^{k+1}})$ by Lemma
\ref{lem-divisibility}, hence $A^2 \equiv I \mod {2^{k+1}}.$  Moreover, $A \in \Gamma(2)$ and
$A = \left( \begin{smallmatrix}
I_q & 0 \\ * & * \end{smallmatrix} \right).$   Let $p \in \mathbf{GL}(n, \mathbb Z)$ be the
change of basis provided by Lemma \ref{lem-GLconj}.  Then $p^{-1}Ap = \left(
\begin{smallmatrix} I_q & 0 \\ * & * \end{smallmatrix} \right)$ and (after re-ordering the
coordinates if necessary), $p^{-1}Ap \equiv \left( \begin{smallmatrix} I_r & 0 \\ 0 & -I_s
\end{smallmatrix} \right) \mod {2^k} $ for some $ r \ge q.$   Set $h = \left(
\begin{smallmatrix} p & 0 \\ 0 & \tr p^{-1} \end{smallmatrix} \right) \in \ker(\nu).$  Set $g
= ah.$  
\begin{equation*} \begin{CD}
F_q @>{\tilde h}>> F_q @>{\tilde a}>> F \\ @AAA @A{u}AA @AA{\gamma'}A \\
F_q @>>{h}> F_q @>>{a}> F \end{CD}. \end{equation*}
Then $\tilde h^{-1} u h = \tilde g^{-1} \gamma ' g \in \ker(\nu)$ and
 $\tilde g^{-1} \gamma' g$ has the desired form (\ref{eqn-normalform}).  \qed

\begin{thm}  Let $(F,\gamma)$ be a real boundary pair. Then there exists $\gamma_1 \in
\Gamma(2^k)$ so that the set $F^{\gamma} = F^{\gamma_1}$ of $\gamma$-real points is
contained in the closure $\overline{\mathfrak h_n^{\gamma_1}}.$  \end{thm}

\subsection{Proof} Set $q = \text{rank}(F).$   By Lemma \ref{lem-IrIs}, there exists $\gamma_1
\in \Gamma(2^k)$ and there exists $g \in \mathbf{Sp}(2n,\mathbb Z)$ so that $F^{\gamma} =
F^{\gamma_1}$, so that $g(F_{n-1}) = F$ and so that $u = \tilde g^{-1} \gamma_1 g$ lies in
$\ker(\nu)$ and has the form (\ref{eqn-normalform}), for some $r \ge q.$  Therefore $j_ruj_r
\equiv I \mod {2^k}$ so Proposition \ref{prop-killer} may be applied.  \qed

\subsection{Proof of Lemma \ref{lem-GLconj}}\label{subsec-proofoflemma}
The lemma is equivalent to the following statement.  Suppose $M$ is a free $\mathbb Z$-module
of rank $n$.  Let $\alpha : M \to M$ be an automorphism such that $(\alpha -I)M \subset 2M$
and $(\alpha^2 -I)M \subset {2^{k+1}}M,$ that is, $\alpha \equiv I \mod 2$ and $\alpha^2
\equiv I \mod {2^{k+1}}.$  Then there exists a basis $x_1,x_2, \ldots, x_n$ of $M$ so
that $\alpha(x_i)=\pm x_i \in {2^k}M$ for $i = 1,2,\ldots, n.$

This statement will be proven by induction on the rank of $M.$  The case of rank 1 is obvious,
so suppose that $M$ has rank $n.$  We will show that there exists a basis $\{x_1,x_2, \ldots,
x_n\}$ of $M$ so that $\alpha(x_1) \equiv \pm x_1 \mod {2^k}.$  If $\ker(\alpha - I)$ is not
trivial or if $\ker(\alpha +I)$ is not trivial then any primitive $x$ within this kernel may
be extended to a basis.  Therefore we may assume that both $(\alpha -I)M$ and $(\alpha +I)M$ 
have maximal rank. 

By elementary divisor theory there exists a basis $x_1,x_2,\ldots,x_n$ of $M$ and integers
$d_1,d_2,$ $\ldots,d_n$ so that $d_1 | d_2 | \ldots | d_n$ and so that $d_1x_1,d_2x_2,\ldots,
d_nx_n$ is a basis of $(\alpha -I)M.$  Similarly there exists another basis
$y_1,y_2,\ldots,y_n$ of $M$ and integers $e_1|e_2|\ldots |e_n$ so that $e_1y_1, e_2y_2,$
$\ldots, e_ny_n$ is a basis for the submodule $(\alpha +I)M\subset M.$

We claim that either $\alpha(x_1) \equiv - x_1 \mod {2^k}$ or $\alpha(y_1) \equiv  y_1 \mod
{2^k}.$  First note that $(\alpha +I)d_1x_1 \equiv 0 \mod {2^{k+1}}$ since $d_1x_1 \in (\alpha
-I)M.$   If $d_1$ is odd this implies $\alpha x_1 + x_1 \equiv 0 \mod {2^{k+1}}$.  If $d_1/2$
is odd, it implies that $\alpha x_1 + x_1 \equiv 0 \mod {2^k}.$  Similarly, if $e_1$ is odd or
if $e_1/2$ is odd then $\alpha y_1 - y_1 \equiv 0 \mod {2^k}.$  However,
 $(\alpha -I)M + (\alpha +I)M = 2M$ so the highest power of 2 which divides 
$\gcd(d_1,e_1)$ is $2^1.$  Therefore one of these four cases must occur, which proves the
claim.

By switching the $x$'s with the $y$'s if necessary, we arrive at a basis $x_1,x_2, \ldots,
x_n$ of $M$ such that $\alpha(x_1) \equiv \pm x_1 \mod {2^k}$.  Write $M = M_1 \oplus M_2$
where $M_1 =\mathbb Zx_1$ and $M_2 = \sum_{i\ge 2} \mathbb Z x_i.$  With respect to this
decomposition, $\alpha$ has the matrix $\left( \begin{smallmatrix} \alpha_{11} & \alpha_{12}\\
\alpha_{21} & \alpha_{22} \end{smallmatrix} \right)$ where $\alpha_{11} \equiv \pm 1 \mod
{2^k}$, and where $\alpha_{22}:M_2 \to M_2.$  We
claim that $\alpha_{22} \equiv I \mod 2$ and $\alpha_{22}^2 \equiv I \mod {2^{k+1}}.$    Since
$\alpha \equiv I \mod 2$  we have
\begin{equation*}
\alpha \left(\begin{matrix} 0 \\ m_2 \end{matrix}\right) = \left(
\begin{matrix} \alpha_{12}m_2 \\ \alpha_{22} m_2 \end{matrix} \right) \equiv
\left( \begin{matrix} 0 \\ m_2 \end{matrix} \right) \mod 2 \end{equation*}
hence $\alpha_{12} \equiv 0 \mod 2$ and $\alpha_{22} \equiv I \mod 2.$  Also,
\begin{equation*}
\alpha\left(\begin{matrix}x_1 \\ 0 \end{matrix} \right) = \left(
\begin{matrix} \alpha_{11}x_1 \\ \alpha_{21}x_1 \end{matrix} \right)
\equiv \left( \begin{matrix} \pm x_1 \\ 0 \end{matrix} \right) \mod {2^k}
\end{equation*}
hence $\alpha_{21} \equiv 0 \mod {2^k}.$
Similarly,
\begin{equation*}
\alpha^2 \left( \begin{matrix} 0 \\ m_2 \end{matrix} \right) = 
\left( \begin{matrix} * \\ \alpha_{21}\alpha_{12}m_2 + \alpha_{22}^2 m_2
\end{matrix}
\right) \equiv \left( \begin{matrix} 0 \\ m_2 \end{matrix} \right) \mod {2^{k+1}}. 
\end{equation*}
But $\alpha_{12} \equiv 0 \mod 2$ and $\alpha_{21} \equiv 0 \mod {2^k}$ so the first term in
this sum is congruent to $0 \mod {2^{k+1}},$ hence $\alpha_{22}^2m_2 \equiv m_2 \mod
{2^{k+1}}$ as claimed.

Therefore we may apply induction to the pair $(M_2, \alpha_{22})$ to obtain a basis, which we
again denote by $x_2,x_3,\ldots, x_n$ such that $\alpha_{22}(x_j) \equiv \pm x_j \mod {2^k}$
for $j \ge 2.$    Hence, for $j \ge 2,$ there are integers $a_j$ and a sign $\epsilon_j = \pm
1$ so that $\alpha(x_j) \equiv a_jx_1 + \epsilon_jx_j \mod {2^k}.$  Each $a_j$ is even since
$\alpha \equiv I
\mod 2.$  Define a new basis
\begin{equation*}
x'_j = \begin{cases} x_j & \text{if } \epsilon_j = \epsilon_1 \\
     \epsilon_j x_j + \frac{1}{2}a_jx_1 &\text{if } \epsilon_j = -\epsilon_1.
\end{cases} 
\end{equation*} 
We claim this basis has the desired property:  $\alpha(x'_j) = \pm x'_j \mod {2^k}.$  First
suppose $\epsilon_j = +\epsilon_1.$  Then $\alpha^2(x_j) = 2a_j\epsilon_jx_1 + x_j.$  Since
$\alpha^2 \equiv I \mod {2^{k+1}}$ we see that $a_k \equiv 0 \mod {2^k}.$  Hence $\alpha(x'_j)
\equiv x_j \mod {2^k}$ as desired.  If $\epsilon_j \epsilon_1 = -1$ then 
\begin{align*}
\alpha(x'_j) &=\epsilon_j\alpha(x_j) + \frac{1}{2} a_j \alpha(x_1) \\
&\equiv \epsilon_j(a_jx_1 + \epsilon_jx_j) + \frac{1}{2}a_j\epsilon_1x_1 \mod{2^k} \\
&\equiv x_j + \frac{1}{2}\epsilon_ja_jx_1 \mod{2^k} \\
&= \epsilon_j(\epsilon_jx_j +\frac{1}{2}a_jx_1) = \epsilon_jx'_j.\end{align*}
  This completes the construction of the desired basis.

To prove the ``moreover'' part of the lemma, let $M' = \sum_{i=1}^q \mathbb Z e_i$ be the
submodule of $M$ generated by the first $q$ standard basis vectors.  Apply the lemma to the
quotient module $\alpha: M/M' \to M/M'.$  Choose any lift $x_{q+1}, \ldots, x_n$ of the
resulting basis of $M/M'$ to $M$ and define $x_1=e_1, x_2=e_2,\ldots, x_q=e_q.$  With respect
to this basis, 
\begin{equation*}
\alpha(x_j) = \epsilon_jx_j + \sum_{j=1}^q a_{ij}x_j \text{ for }i > q\end{equation*}
where $\epsilon_i = \pm 1$ and $a_{ij}$ is even.  Set 
\begin{equation*}
x'_i = \begin{cases} x_i &\text{if } \epsilon_i = 1\\
-x_i + \sum_{j=1}^q \frac{1}{2}a_{ij}x_j &\text{if } \epsilon_i = -1
\end{cases}. \end{equation*}
Then $\alpha(x'_i) \equiv \epsilon_ix'_i \mod {2^k}$ and the change of basis matrix has the
desired form.  \qed

\end{document}